\documentclass[11pt,a4paper]{amsart}
\usepackage[a4paper,margin=1in]{geometry}
\usepackage[T1]{fontenc}
\usepackage[utf8]{inputenc}
\usepackage{lmodern}
\usepackage{amsmath,amssymb,amsthm,mathtools,bm}
\usepackage{booktabs}
\usepackage{graphicx}
\usepackage{mathrsfs}
\usepackage{enumitem}
\usepackage[protrusion=true,expansion=false]{microtype}
\usepackage[hidelinks]{hyperref}

\theoremstyle{plain}
\newtheorem{theorem}{Theorem}[section]
\newtheorem{lemma}[theorem]{Lemma}
\newtheorem{proposition}[theorem]{Proposition}

\newtheorem{assumption}[theorem]{Assumption}
\theoremstyle{definition}
\newtheorem{definition}[theorem]{Definition}
\newtheorem{algorithm}[theorem]{Algorithm}
\theoremstyle{remark}
\newtheorem{remark}[theorem]{Remark}

\newcommand{\Ga}{\Gamma}
\newcommand{\D}{\mathcal D}
\newcommand{\A}{\mathcal A}
\newcommand{\C}{\mathbb C}
\newcommand{\R}{\mathbb R}

\newcommand{\spanop}{\operatorname{span}}
\newcommand{\dist}{\operatorname{dist}}

\newcommand{\norm}[1]{\left\|#1\right\|}

\title[B-spline--Heaviside Collocation]{Periodic B-spline--Heaviside Collocation for Fredholm Integral Equations with Piecewise H\"older Data}

\author{Maria Capcelea}
\address{Institute of Mathematics and Computer Science, Moldova State University, Chisinau, Republic of Moldova}
\email{maria.capcelea@usm.md}

\author{Titu Capcelea}
\address{Department of Computer Science, Moldova State University, Chisinau, Republic of Moldova}
\email{titu.capcelea@usm.md}

\date{}

\begin{document}

\begin{abstract}
We propose and analyze a periodic B-spline--Heaviside collocation method for Fredholm integral equations of the second kind on a smooth closed contour with piecewise H\"older data and finitely many prescribed jumps. Since the regular integral term is globally continuous, the solution jumps coincide with those of the data. The approximation space combines one global periodic B-spline component with cut-compensated Heaviside generators carrying exactly the prescribed physical jumps, which yields a triangular algorithm: the jump amplitudes are reconstructed first, and only the continuous component is obtained from an $n_B\times n_B$ collocation system. For B-spline orders $m=2,3,4$, exact collocation is uniformly stable and converges in $PH_\beta$, $0<\beta<\alpha<\mu\le1$, with order $O(h^{\alpha-\beta})$. If the parametrized kernel density is H\"older continuous of exponent $\rho$ in the integration variable, a panelwise $q$-point Gauss rule with $q\ge\lceil m/2\rceil$ produces an $O(h^\rho)$ quadrature perturbation on the discrete spline space and preserves the exact-collocation rate whenever $\rho\ge\alpha-\beta$. The exact-operator one-step B-spline--Heaviside extension satisfies the stronger estimate $O(h^\alpha)$ in $H_\mu$. We also treat sampled lateral data and a nonvanishing piecewise H\"older leading coefficient. Numerical experiments are consistent with the predicted stability and convergence, show the substantial error reduction produced by iteration, and indicate competitive accuracy relative to established discontinuity-adapted collocation and Nystr\"om alternatives, while the proposed representation introduces no mesh-dependent artificial jumps.

\end{abstract}

\maketitle

\section{Introduction}\label{sec:introduction}

Fredholm integral equations of the second kind arise naturally in boundary integral formulations, potential theory, elasticity, inverse problems, and transport models.  Their numerical solution by projection, quadrature, and collocation methods has a long history; standard references include \cite{Atkinson,Kress}.  Collocation remains an active area of research, especially when the geometry is noncompact or the kernel and solution have reduced regularity.  Recent developments include multistep schemes, mapped spectral collocation on infinite intervals, radial-basis-function discretizations, weighted Jacobi constructions, and mapped Legendre methods \cite{Wang2022,Rahmoune2021,Molabahrami2023,Allouch2024,Boutarcha2026}.

A particular difficulty occurs when the data or the solution is not globally smooth.  Piecewise-continuous collocation and its iterated form were developed by Atkinson, Graham, and Sloan \cite{AtkinsonGrahamSloan}; Joe studied discontinuous piecewise-polynomial and fully discrete collocation \cite{Joe1985,JoeDiscrete1985}, and iterated spline collocation was further analyzed in \cite{GrahamJoeSloan1985}.  Related approximation strategies for nonsmooth functions include nonlinear B-spline-type approximants, WENO modifications of spline quasi-interpolants, Fourier-informed knot placement, and correction terms for known singularities \cite{Amat2023,Arandiga2023,Lenz2023,LiPanRuizYanez2024}.  Within this line of research, the present authors proposed a B-spline--Heaviside enrichment for the approximation of discontinuous functions on closed contours and for Fredholm collocation with discontinuous data \cite{CapceleaCapcelea2022,CapceleaCapcelea2023}.  The present paper develops this approach into a periodic collocation framework with a complete stability, convergence, quadrature, and iterated-extension analysis.

Against this background, let $\Ga\subset\C$ be a smooth simple closed contour and consider
\begin{equation}\label{eq:main}
        (A\varphi)(t)
        :=\varphi(t)-\lambda(K\varphi)(t)=f(t),
        \qquad
        (K\varphi)(t):=\int_\Ga k(t,\tau)\varphi(\tau)\,d\tau.
\end{equation}
We study the case in which the right-hand side and the solution are piecewise H\"older continuous and may have finitely many prescribed jumps on $\Ga$.  Let
$
        \D=\{t_1^d,\ldots,t_{n_d}^d\}\subset\Ga
$
denote the corresponding finite set of distinct physical jump points, ordered according to the positive orientation of the contour.  Two features make this problem different from standard smooth periodic collocation.  First, a globally continuous periodic spline space cannot reproduce a nonzero jump uniformly.  Second, a cut-based Heaviside function on a closed contour carries not only the intended physical jump but also a compensating discontinuity at the parameter cut.

The central idea of the paper is to use the enriched periodic space
\[
        V_h^H
        =
        S_h\oplus
        \spanop\{G_1,\ldots,G_{n_d}\},
        \qquad
        G_j=H_j-R,
\]
where $S_h$ is one global periodic B-spline space and $R$ is a parameter ramp.  Each compensated generator $G_j$ has a unit jump at the prescribed point $t_j^d$ and no artificial jump at the cut.  The enriched space is the core approximation mechanism of the method: the B-spline component represents the globally continuous behavior, while the corrected Heaviside generators represent exactly the prescribed discontinuities.

The construction is reinforced by a structural property of regular Fredholm operators.  Since $K\varphi$ is globally continuous,
\[
        [\varphi]_{t_j^d}=[f]_{t_j^d},
        \qquad t_j^d\in\D.
\]
Thus the jump amplitudes are determined directly by the data and need not be introduced as additional unknowns in the dense Fredholm system.  The resulting algorithm is triangular: the jump component is reconstructed from exact or estimated one-sided values, removed from the equation, and then one periodic B-spline collocation problem is solved for the continuous component.  Finally, the prescribed jumps are restored.  This reduction preserves the global periodic spline representation while separating the discontinuous information before the dense Fredholm system is solved.

Within this framework, the results developed in the paper may be summarized as follows.
\begin{enumerate}[label=(\roman*),leftmargin=1.1cm]
\item We construct cut-compensated periodic jump generators and prove the corresponding unique continuous--jump decomposition.
\item We derive the exact triangular reduction and a periodic B-spline collocation method of orders $m=2,3,4$ whose dense linear system has dimension $n_B$, independently of the number of prescribed jumps.
\item Under
\[
        0<\beta<\alpha<\mu\le1,
\]
a uniform H\"older smoothing condition on the kernel, and a nonresonance condition on $\lambda$, we establish uniform discrete stability and
\[
        \norm{\varphi-\varphi_h}_{PH_\beta}
        \le
        C_\lambda h^{\alpha-\beta}
        \norm{f}_{PH_\alpha}.
\]
\item We obtain an explicit fully discrete estimate.  If the parametrized kernel density is $\rho$-H\"older in the integration variable with values in $H_\mu(\Ga)$, then a panelwise Gauss rule with $q\ge\lceil m/2\rceil$ contributes $O(h^\rho)$ and preserves the exact-collocation rate whenever $\rho\ge\alpha-\beta$.
\item For the one-step extension
\[
        \varphi_h^{\mathrm{it}}
        =f+\lambda K\varphi_h,
\]
we prove the improved estimate
\[
        \norm{\varphi-\varphi_h^{\mathrm{it}}}_{H_\mu}
        \le
        C_\lambda h^\alpha\norm{f}_{PH_\alpha}.
\]
\end{enumerate}
The analysis is complemented by one-sided reconstruction from sampled data, an extension to a nonvanishing piecewise H\"older leading coefficient, and a reproducible numerical comparison with established discontinuity-adapted collocation and Nystr\"om alternatives.

The paper is organized as follows.  Section~\ref{sec:spaces} defines the contour geometry and the piecewise H\"older spaces.  Section~\ref{sec:operator} establishes the mapping properties of the Fredholm operator, solvability, and the jump identity.  Section~\ref{sec:lateral} treats the reconstruction of lateral data from samples.  Sections~\ref{sec:generators} and~\ref{sec:reduction} construct the corrected periodic generators and derive the reduced continuous equation.  Section~\ref{sec:splines} introduces the periodic B-spline spaces, and Section~\ref{sec:convergence} proves stability and convergence of exact collocation.  Sections~\ref{sec:iterated} and~\ref{sec:quadrature} develop the iterated extension and the fully discrete quadrature analysis.  Section~\ref{sec:variablecoefficient} treats a piecewise H\"older leading coefficient, Section~\ref{sec:algorithm} summarizes the implementation algorithm, Section~\ref{sec:numerics} presents the numerical experiments, and Section~\ref{sec:conclusions} gives the conclusions.

\section{Contour geometry and piecewise H\"older spaces}\label{sec:spaces}

Let $\Ga\subset\C$ be a simple, closed, positively oriented contour of class $C^2$.  Fix a regular $2\pi$-periodic parametrization
\[
        \gamma:\R\to\Ga,
        \qquad \gamma(\theta+2\pi)=\gamma(\theta),
        \qquad \gamma'(\theta)\ne0.
\]
Let
$
        \D=\{t_1^d,\ldots,t_{n_d}^d\}\subset\Ga
$
be distinct points ordered according to the positive orientation.  Choose an auxiliary cut point $t_*\in\Ga\setminus\D$ and shift the parameter so that
$
        \gamma(0)=\gamma(2\pi)=t_*.
$
Then
\[
        t_j^d=\gamma(\theta_j^d),
        \qquad 0<\theta_1^d<\cdots<\theta_{n_d}^d<2\pi.
\]
The open arcs between consecutive jump points are denoted by
$
        \A(\Ga\setminus\D)=\{\Ga_1,\ldots,\Ga_{n_d}\},
$
with cyclic indexing.  For a function with finite one-sided limits, set
$
        [v]_{t_j^d}:=v(t_j^d+0)-v(t_j^d-0).
$
We represent functions by their incoming values,
$
        v(t_j^d)=v(t_j^d-0).
$

\begin{definition}\label{def:Xsigma}
For $0<\sigma<1$, the space
\[
        X_\sigma:=PH_\sigma(\Ga,\D)
\]
consists of functions that are $\sigma$-H\"older continuous on every continuity arc and possess finite lateral limits at all points of $\D$.  For an arc $I\in\A(\Ga\setminus\D)$, let $\overline I^{\,\ell}$ denote the arc completed by the appropriate lateral endpoint values, and define
\[
        \norm{v}_{\infty,I}
        :=\sup_{t\in\overline I^{\,\ell}}|v(t)|,
        \qquad
        [v]_{\sigma,I}
        :=\sup_{\substack{t,s\in\overline I^{\,\ell}\\t\ne s}}
        \frac{|v(t)-v(s)|}{|t-s|^\sigma}.
\]
The norm in $X_\sigma$ is
\begin{equation}\label{eq:Xnorm}
        \norm{v}_{X_\sigma}
        :=\max_{I\in\A(\Ga\setminus\D)}
        \left(\norm{v}_{\infty,I}+[v]_{\sigma,I}\right).
\end{equation}
The subspace of globally $\sigma$-H\"older continuous functions is denoted by $H_\sigma(\Ga)$ and is equipped with the analogous global norm
$
        \norm{v}_{H_\sigma(\Ga)}
        :=\norm{v}_{\infty,\Ga}+[v]_{\sigma,\Ga}.
$
\end{definition}

For $\sigma=1$, $H_1(\Ga)$ denotes the Lipschitz space equipped with the analogous norm.

The norm \eqref{eq:Xnorm} records the supremum and H\"older seminorm separately on every continuity arc.  The Euclidean, arclength, and parameter metrics are locally equivalent because $\gamma$ is regular and $C^2$.  The corresponding H\"older norms are therefore equivalent.  Both $X_\sigma$ and $H_\sigma(\Ga)$ are Banach spaces.  If $0<\beta<\alpha<1$, then
\begin{equation}\label{eq:embedding}
        X_\alpha\hookrightarrow X_\beta,
        \qquad H_\alpha(\Ga)\hookrightarrow H_\beta(\Ga).
\end{equation}

The embeddings in \eqref{eq:embedding} will be used repeatedly.  We impose the following smoothing assumption on the regular kernel.

\begin{assumption}[Kernel regularity]\label{ass:kernel}
There exists an exponent $\mu$ such that
\[
        0<\beta<\alpha<\mu\le1
\]
and constants $M_0,M_\mu$ for which
\begin{align}
        |k(t,\tau)|&\le M_0,\label{eq:kbound}\\
        |k(t,\tau)-k(s,\tau)|&\le M_\mu|t-s|^\mu\label{eq:kholder}
\end{align}
for all $t,s,\tau\in\Ga$.  In addition, the family $\{\tau\mapsto k(t,\tau):t\in\Ga\}$ is equicontinuous.  Equivalently, there is a modulus of continuity $\omega_k$, with $\omega_k(r)\to0$ as $r\downarrow0$, such that
\begin{equation}\label{eq:kuniformtau}
        |k(t,\tau)-k(t,\zeta)|
        \le \omega_k(|\tau-\zeta|),
        \qquad t,\tau,\zeta\in\Ga.
\end{equation}
\end{assumption}

The bound \eqref{eq:kbound} and the estimate \eqref{eq:kholder} are the assumptions used in the mapping, compactness, and exact-collocation arguments below.  Condition \eqref{eq:kholder} is uniform H\"older continuity in the observation variable $t$, uniformly with respect to $\tau$, whereas \eqref{eq:kuniformtau} is uniform continuity in the integration variable $\tau$, uniformly with respect to $t$; the latter is included to justify standard panelwise quadrature for the original kernel.  The condition $\mu>\alpha$ is stronger than is needed to define the equation, but it gives a clean operator-norm convergence theory in the H\"older spaces used below.  A jointly continuous kernel that is $C^1$ in its first variable satisfies Assumption~\ref{ass:kernel} for every $\mu\le1$.

\section{Fredholm mapping properties and the jump identity}\label{sec:operator}

\begin{lemma}[Smoothing and compactness]\label{lem:Kmapping}
Under Assumption~\ref{ass:kernel},
\[
        K:X_\beta\longrightarrow H_\mu(\Ga)
\]
is bounded.  Consequently, $K:X_\beta\to X_\beta$ is compact.
\end{lemma}

\begin{proof}
Let $L_\Ga$ be the length of the contour.  For $v\in X_\beta$,
$
        |(Kv)(t)|
        \le L_\Ga M_0\norm{v}_{\infty,\Ga},
$
and
$
        |(Kv)(t)-(Kv)(s)|
        \le L_\Ga M_\mu\norm{v}_{\infty,\Ga}|t-s|^\mu.
$
Thus
$
        \norm{Kv}_{H_\mu(\Ga)}
        \le C\norm{v}_{X_\beta}.
$
The embedding $H_\mu(\Ga)\hookrightarrow H_\beta(\Ga)$ is compact for $\beta<\mu$.  Since $H_\beta(\Ga)$ is continuously embedded in $X_\beta$, the asserted compactness follows.
\end{proof}

\begin{assumption}[Nonresonance of the parameter]\label{ass:injectivity}
The fixed parameter $\lambda\in\C$ is such that
\begin{equation}\label{eq:nonresonance}
        \ker(I-\lambda K)=\{0\}
        \qquad\text{in }X_\beta.
\end{equation}
Equivalently, either $\lambda=0$, or
\[
        \lambda^{-1}\notin\sigma(K:X_\beta\to X_\beta).
\]
\end{assumption}

\begin{remark}[Admissible values of $\lambda$]\label{rem:lambda}
By Lemma~\ref{lem:Kmapping}, $K:X_\beta\to X_\beta$ is compact.  Hence every nonzero point of its spectrum is an eigenvalue, and the nonresonance condition \eqref{eq:nonresonance} excludes only the characteristic values
$
        \lambda=\nu^{-1},
         0\ne\nu\in\sigma(K).
$
These excluded values are isolated and can accumulate only at infinity.  A simple sufficient, but not necessary, condition is
\begin{equation}\label{eq:neumannlambda}
        |\lambda|\,\norm{K}_{X_\beta\to X_\beta}<1,
\end{equation}
Condition \eqref{eq:neumannlambda} implies invertibility by the Neumann series.  No universal restriction such as $|\lambda|<1$ is required.  The relevant stability constants may grow when $\lambda$ approaches a characteristic value.
\end{remark}

\begin{theorem}[Well-posedness and regularity]\label{thm:wellposed}
Under Assumptions~\ref{ass:kernel} and~\ref{ass:injectivity}, the operator
\[
        A=I-\lambda K:X_\beta\to X_\beta
\]
is boundedly invertible.  If $f\in X_\alpha$, then the unique solution of \eqref{eq:main} belongs to $X_\alpha$ and
\begin{equation}\label{eq:continuous-stability}
        \norm{\varphi}_{X_\alpha}
        \le C_\lambda\norm{f}_{X_\alpha}.
\end{equation}
Here $C_\lambda$ may depend on the contour, the kernel, the H\"older exponents, and $\lambda$, in particular through the resolvent norm
$\norm{(I-\lambda K)^{-1}}_{X_\beta\to X_\beta}$, but it is independent of $f$.
\end{theorem}

\begin{proof}
By Lemma~\ref{lem:Kmapping}, $K$ is compact on $X_\beta$.  Hence $A=I-\lambda K$ is a Fredholm operator of index zero.  Assumption~\ref{ass:injectivity} and the Fredholm alternative imply bounded invertibility on $X_\beta$.

For $f\in X_\alpha$, the solution initially satisfies
\[
        \norm{\varphi}_{X_\beta}
        \le C_\lambda\norm{f}_{X_\beta}
        \le C_\lambda\norm{f}_{X_\alpha}.
\]
Equation \eqref{eq:main} gives
$
        \varphi=f+\lambda K\varphi.
$
The first term belongs to $X_\alpha$, while Lemma~\ref{lem:Kmapping} gives $K\varphi\in H_\mu(\Ga)\subset H_\alpha(\Ga)$.  Thus $\varphi\in X_\alpha$, and the same estimates yield \eqref{eq:continuous-stability}.
\end{proof}

The well-posedness estimate \eqref{eq:continuous-stability} will also control data perturbations.  The decisive simplification of the Fredholm problem is the following identity.

\begin{proposition}[Exact jump transfer]\label{prop:jumpidentity}
For every $v\in X_\beta$ and every $t_j^d\in\D$,
\begin{equation}\label{eq:jumpidentity}
        [Av]_{t_j^d}=[v]_{t_j^d}.
\end{equation}
Consequently, the solution of \eqref{eq:main} satisfies
\begin{equation}\label{eq:solutionjumps}
        [\varphi]_{t_j^d}=[f]_{t_j^d},
        \qquad j=1,\ldots,n_d.
\end{equation}
\end{proposition}

\begin{proof}
Lemma~\ref{lem:Kmapping} shows that $Kv$ is globally continuous.  Therefore $[Kv]_{t_j^d}=0$, and \eqref{eq:jumpidentity} follows from $Av=v-\lambda Kv$.  Applying this identity to $A\varphi=f$ proves \eqref{eq:solutionjumps}.
\end{proof}

\begin{remark}[Required input at a jump]\label{rem:jumpdata}
The solution-jump formula \eqref{eq:solutionjumps} requires genuine lateral data: a single assigned value $f(t_j^d)$ does not determine the jump.  The method requires either the two lateral values $f(t_j^d-0)$ and $f(t_j^d+0)$ or their difference.  When only sampled values are available, the two traces must be reconstructed separately, using points from the two adjacent continuity arcs and never mixing observations across the jump.
\end{remark}

\section{Reconstruction of lateral data from samples}\label{sec:lateral}

Assume that the jump parameters $\theta_j^d$ are known and that the available data are
\begin{equation}\label{eq:samples}
        y_\ell=f(\gamma(\vartheta_\ell))+e_\ell,
        \qquad \ell=1,\ldots,N,
\end{equation}
where the sample parameters $\vartheta_\ell$ are ordered cyclically, none of them belongs to $\{\theta_j^d\}$, and $e_\ell$ represents measurement or rounding error.  For every jump, choose one-sided bandwidths $r_j^-,r_j^+>0$ small enough that the windows
\[
 (\theta_j^d-r_j^-,\theta_j^d),
 \qquad
 (\theta_j^d,\theta_j^d+r_j^+)
\]
remain inside the two adjacent continuity arcs.  All parameter differences in this section are understood cyclically after the cut has been placed away from $\D$.

Let $W:[0,1]\to[0,\infty)$ be a bounded weight function, positive on a neighborhood of the origin.  Define the index sets
\[
 \mathcal I_j^-:=\{\ell:0<\theta_j^d-\vartheta_\ell\le r_j^-\}, \quad
 \mathcal I_j^+:=\{\ell:0<\vartheta_\ell-\theta_j^d\le r_j^+\},
\]
and the normalized one-sided weights
\[
 w_{j\ell}^-:=
 \frac{W((\theta_j^d-\vartheta_\ell)/r_j^-)}
 {\sum_{q\in\mathcal I_j^-}W((\theta_j^d-\vartheta_q)/r_j^-)},\quad
 w_{j\ell}^+:=
 \frac{W((\vartheta_\ell-\theta_j^d)/r_j^+)}
 {\sum_{q\in\mathcal I_j^+}W((\vartheta_q-\theta_j^d)/r_j^+)}.
\]
Provided that both denominators are nonzero, the one-sided local-constant estimators are
\begin{equation}\label{eq:lateralestimators}
 \widehat f_j^-:=\sum_{\ell\in\mathcal I_j^-}w_{j\ell}^-y_\ell,
 \qquad
 \widehat f_j^+:=\sum_{\ell\in\mathcal I_j^+}w_{j\ell}^+y_\ell,
 \qquad
 \widehat\delta_j:=\widehat f_j^+-\widehat f_j^-.
\end{equation}
For noise-free data, the nearest sample on each side is the simplest admissible choice.  In the presence of noise, several samples should be averaged, while the bandwidths must still be small enough to control the one-sided extrapolation bias.

\begin{proposition}[Accuracy of lateral reconstruction]\label{prop:lateralerror}
Suppose that $f\in X_\alpha$ and $|e_\ell|\le\varepsilon$.  Then the estimators in \eqref{eq:lateralestimators} satisfy
\begin{align}
 |\widehat f_j^- - f(t_j^d-0)|
 &\le C_\gamma[f]_{\alpha,\Ga_j^-}(r_j^-)^\alpha+\varepsilon,\label{eq:lefttraceerror}\\
 |\widehat f_j^+ - f(t_j^d+0)|
 &\le C_\gamma[f]_{\alpha,\Ga_j^+}(r_j^+)^\alpha+\varepsilon,\label{eq:righttraceerror}\\
 |\widehat\delta_j-[f]_{t_j^d}|
 &\le C_\gamma\left([f]_{\alpha,\Ga_j^-}(r_j^-)^\alpha
      +[f]_{\alpha,\Ga_j^+}(r_j^+)^\alpha\right)+2\varepsilon,
 \label{eq:jumperror}
\end{align}
where $\Ga_j^-$ and $\Ga_j^+$ are the adjacent incoming and outgoing arcs, and $C_\gamma$ depends only on the equivalence between the parameter and contour metrics.
\end{proposition}

\begin{proof}
The weights on each side are nonnegative and sum to one.  Every sample used in the left estimator is at parameter distance at most $r_j^-$ from the left trace.  The equivalence of the parameter and contour metrics, the one-sided H\"older estimate, and the error bound in \eqref{eq:samples} give \eqref{eq:lefttraceerror}.  The right estimate \eqref{eq:righttraceerror} is identical, and \eqref{eq:jumperror} follows by subtraction.
\end{proof}

\begin{remark}[Consistency of the jump estimator]\label{rem:lateral-consistency}
The right-hand side of \eqref{eq:jumperror} tends to zero in a sampling regime for which $r_j^-\to0$, $r_j^+\to0$, and the error level $\varepsilon=\varepsilon_N\to0$, provided that each shrinking one-sided window still contains sufficiently many samples.  For a fixed nonzero noise bound $\varepsilon$, the deterministic bias tends to zero but the estimate has the irreducible upper bound $2\varepsilon$.  Thus consistency requires both increasing one-sided sampling density and vanishing measurement or rounding error (or an averaging procedure with a vanishing effective noise level).
\end{remark}

\begin{algorithm}[Sample-based reconstruction of the jump data]\label{alg:lateral}
For every $t_j^d\in\D$, perform the following steps.
\begin{enumerate}[label=\textbf{Step \arabic*.},leftmargin=2.6cm]
\item Sort the samples by the periodic parameter and identify the two continuity arcs adjacent to $\theta_j^d$.
\item Select either a fixed number $q_-$ and $q_+$ of nearest samples on the two sides or bandwidths $r_j^-$ and $r_j^+$.  Reject any window that crosses another jump or the parameter cut.
\item Compute $\widehat f_j^-$ and $\widehat f_j^+$ from \eqref{eq:lateralestimators}.  Record the largest one-sided sample distances, since they control the deterministic bias in Proposition~\ref{prop:lateralerror}.
\item Set $\widehat\delta_j=\widehat f_j^+-\widehat f_j^-$.  Repeat the computation with a smaller window or fewer nearest points and use the change in $\widehat\delta_j$ as a practical reconstruction diagnostic.
\item If the collocation values $f(t_i^B)$ are not sampled directly, reconstruct them by interpolation or local regression using only samples from the continuity arc containing $t_i^B$.  Once the corrected generators of Section~\ref{sec:generators} have been constructed, form
\[
        \widehat f_C(t_i^B)
        :=\widehat f(t_i^B)-\sum_{j=1}^{n_d}\widehat\delta_jG_j(t_i^B).
\]
\end{enumerate}
\end{algorithm}

For data known only to be piecewise H\"older, the local-constant construction gives the natural worst-case rate $O(r^\alpha)$; a higher-degree fit cannot improve that guaranteed rate without additional smoothness.  If each one-sided restriction is known to be smoother, a one-sided local polynomial of degree $p$ may be fitted in the scaled distance from the jump and evaluated at the endpoint.  Such a replacement can reduce bias, but it requires at least $p+1$ well-distributed samples on each side and an explicit condition-number check for the weighted least-squares matrix.

The effect of reconstruction error is controlled by the stability of the Fredholm equation.  In particular, if a reconstructed datum $\widehat f\in X_\beta$ satisfies $\norm{\widehat f-f}_{X_\beta}\le\xi$, then the corresponding exact solutions satisfy
\begin{equation}\label{eq:data-perturbation}
        \norm{\widehat\varphi-\varphi}_{X_\beta}
        \le \norm{A^{-1}}_{X_\beta\to X_\beta}\,\xi.
\end{equation}
Estimate \eqref{eq:data-perturbation} shows that the sampling bandwidth and the quadrature tolerance should be chosen so that their errors do not dominate the target discretization scale $h^{\alpha-\beta}$.

\section{Correct periodic jump generators}\label{sec:generators}

For $0\le\theta<2\pi$, define the cut-based one-sided step
\begin{equation}\label{eq:Hj}
 H_j(\gamma(\theta))=
 \begin{cases}
 0,&0\le\theta\le\theta_j^d,\\
 1,&\theta_j^d<\theta<2\pi,
 \end{cases}
 \qquad H_j(t_j^d)=0,
\end{equation}
and the normalized parameter ramp
\begin{equation}\label{eq:ramp}
        R(\gamma(\theta))=\frac{\theta}{2\pi},
        \qquad 0\le\theta<2\pi.
\end{equation}
Both $H_j$ and $R$ have jump $-1$ at the cut $t_*$.  Their difference
\begin{equation}\label{eq:Gj}
        G_j:=H_j-R
\end{equation}
is single-valued and continuous at the cut.

\begin{lemma}[Independent jump coordinates]\label{lem:generatorjumps}
For $i,j=1,\ldots,n_d$,
\begin{equation}\label{eq:jumpmatrix}
        [G_j]_{t_i^d}=\delta_{ij}.
\end{equation}
Moreover, $G_j$ is Lipschitz continuous on every closed subarc not containing $t_j^d$ and belongs to $X_\sigma$ for every $0<\sigma<1$.
\end{lemma}

\begin{proof}
At $t_j^d$, the step $H_j$ has jump $+1$ and $R$ is continuous, which gives $[G_j]_{t_j^d}=1$.  At the other prescribed points both terms are continuous.  At the cut, the incoming and outgoing limits of both $H_j$ and $R$ are respectively $1$ and $0$, so their jumps cancel.  Away from $t_j^d$, the step is constant and $R\circ\gamma$ is smooth in the parameter.
\end{proof}

\begin{remark}[Why the uncorrected step is invalid]\label{rem:rawstep}
The function $H_j$ in \eqref{eq:Hj} cannot be described as a single-jump function on the closed contour: it also jumps at $t_*$.  Therefore a sum $\sum_j\delta_jH_j$ creates an artificial cut jump $-\sum_j\delta_j$.  Unless that sum happens to vanish, subtracting such a combination from a piecewise H\"older function does not produce a continuous remainder.  The ramp $R$ in \eqref{eq:ramp} and the correction \eqref{eq:Gj} are indispensable for arbitrary unbalanced jump vectors.
\end{remark}

\begin{proposition}[Continuous--jump decomposition for a fixed cut]\label{prop:decomposition}
For the fixed cut and parametrization chosen above, let $v\in X_\alpha$ and set
$
        \delta_j=[v]_{t_j^d}.
$
Then
\begin{equation}\label{eq:decomposition}
        v=v_C+J_v,
        \qquad
        J_v:=\sum_{j=1}^{n_d}\delta_jG_j,
        \qquad
        v_C\in H_\alpha(\Ga),
\end{equation}
and the representation is unique.  Moreover,
\begin{equation}\label{eq:decompositionbound}
        \norm{v_C}_{H_\alpha(\Ga)}
        +\sum_{j=1}^{n_d}|\delta_j|
        \le C\norm{v}_{X_\alpha}.
\end{equation}
\end{proposition}

\begin{proof}
Define $v_C=v-J_v$ according to \eqref{eq:decomposition}.  By Lemma~\ref{lem:generatorjumps}, every jump of $v_C$ at $\D$ is zero, and $v_C$ is continuous at the parameter cut.  It is $\alpha$-H\"older on each continuity arc.  A continuous function that is $\alpha$-H\"older on finitely many adjacent arcs is globally $\alpha$-H\"older: split a shorter connecting subarc at the finitely many endpoints and use
\[
        \sum_{r=1}^N a_r^\alpha
        \le N^{1-\alpha}\left(\sum_{r=1}^N a_r\right)^\alpha.
\]
The jump bound $|\delta_j|\le2\norm{v}_{\infty}$ and the fixed H\"older norms of the generators give \eqref{eq:decompositionbound}.

For uniqueness, take the jump at $t_i^d$ in a representation consisting of a continuous function plus a combination of the $G_j$.  Equation \eqref{eq:jumpmatrix} gives the coefficient of $G_i$.  All coefficients therefore vanish, and so does the continuous part.
\end{proof}

\section{Reduction to one continuous Fredholm equation}\label{sec:reduction}

Apply Proposition~\ref{prop:decomposition} to the data and use the decomposition \eqref{eq:fdecomp}:
\begin{equation}\label{eq:fdecomp}
        f=f_C+J_f,
        \qquad
        J_f=\sum_{j=1}^{n_d}[f]_{t_j^d}G_j,
        \qquad f_C\in H_\alpha(\Ga).
\end{equation}
Proposition~\ref{prop:jumpidentity} shows that the exact solution has precisely the same jump component.  Write
\begin{equation}\label{eq:phidecomp}
        \varphi=u+J_f.
\end{equation}
Substitution in \eqref{eq:main} gives the globally continuous equation
\begin{equation}\label{eq:reduced}
        Au=g,
        \qquad
        g:=f_C+\lambda KJ_f.
\end{equation}
By \eqref{eq:fdecomp}--\eqref{eq:reduced}, $f_C\in H_\alpha(\Ga)$ and $KJ_f\in H_\mu(\Ga)$, so $g\in H_\alpha(\Ga)$.

\begin{theorem}[Equivalence of the original and reduced problems]\label{thm:reduction}
Under Assumptions~\ref{ass:kernel} and~\ref{ass:injectivity}, equations \eqref{eq:main} and \eqref{eq:reduced} are equivalent under the transformation \eqref{eq:phidecomp}.  The reduced equation has a unique solution $u\in H_\alpha(\Ga)$, and
\begin{equation}\label{eq:ubound}
        \norm{u}_{H_\alpha(\Ga)}
        \le C\norm{f}_{X_\alpha}.
\end{equation}
\end{theorem}

\begin{proof}
If $\varphi$ solves \eqref{eq:main}, then \eqref{eq:solutionjumps} and the uniqueness of Proposition~\ref{prop:decomposition} imply $\varphi=u+J_f$ with $u$ continuous.  Substitution yields \eqref{eq:reduced}.  Conversely, a solution of \eqref{eq:reduced} gives a solution $\varphi=u+J_f$ of \eqref{eq:main}.

The right-hand side $g$ is globally $\alpha$-H\"older.  Theorem~\ref{thm:wellposed} gives a unique $u\in X_\alpha$.  Since $u=g+\lambda Ku$ and both terms on the right are globally continuous, $u$ has no jumps and therefore lies in $H_\alpha(\Ga)$.  Estimate \eqref{eq:ubound} follows from \eqref{eq:continuous-stability} and Proposition~\ref{prop:decomposition}.
\end{proof}

\begin{remark}[Triangular full system]\label{rem:triangular}
If the jump amplitudes and spline coefficients are formally kept in one vector, the discrete equations have the block form
\[
\begin{pmatrix}
 \mathbb A_h & \mathbb H_h\\
 0 & I_{n_d}
\end{pmatrix}
\begin{pmatrix}\bm a\\\bm\delta\end{pmatrix}
=
\begin{pmatrix}\bm f\\\bm{[f]}\end{pmatrix}.
\]
The lower block gives $\bm\delta=\bm{[f]}$ exactly.  Eliminating it produces the reduced $n_B\times n_B$ system developed below.  Solving the larger system provides no additional information.
\end{remark}

\section{Periodic B-splines of orders two, three, and four}\label{sec:splines}

Fix an order
$
        m\in\{2,3,4\},
$
corresponding respectively to piecewise linear, quadratic, and cubic splines.  Let $n_B\ge4$, put $h=2\pi/n_B$, and choose a phase $\rho_h\in[0,h)$ such that the mesh points stay away from $\D$:
\begin{equation}\label{eq:phase}
        \dist(\theta_j^d-\rho_h,h\mathbb Z)\ge\eta h,
        \qquad j=1,\ldots,n_d,
\end{equation}
where $0<\eta<1/(2n_d)$ is fixed.  Such a phase exists because the union of the forbidden phase intervals has length strictly smaller than $h$.  Define
\[
        \theta_i^B=\rho_h+(i-1)h\pmod{2\pi},
        \qquad t_i^B=\gamma(\theta_i^B),
        \qquad i=1,\ldots,n_B.
\]

Let $N_m$ be the centered cardinal B-spline of order $m$, normalized by $\sum_{r\in\mathbb Z}N_m(x-r)=1$.  For the three orders used here,
\begin{align*}
 N_2(x)&=(1-|x|)_+,\\
 N_3(x)&=
 \begin{cases}
 3/4-x^2,&|x|\le1/2,\\
 \tfrac12(3/2-|x|)^2,&1/2<|x|<3/2,\\
 0,&|x|\ge3/2,
 \end{cases}\\
 N_4(x)&=\frac16
 \begin{cases}
 4-6x^2+3|x|^3,&|x|<1,\\
 (2-|x|)^3,&1\le|x|<2,\\
 0,&|x|\ge2.
 \end{cases}
\end{align*}
The periodic basis functions transported to $\Ga$ are
\begin{equation}\label{eq:periodicB}
        B_{m,j}(\gamma(\theta))
        :=\sum_{\ell\in\mathbb Z}
        N_m\!\left(\frac{\theta-\theta_j^B}{h}-\ell n_B\right),
        \qquad j=1,\ldots,n_B.
\end{equation}
Only finitely many terms are nonzero.  We fix $m$ and write $B_j:=B_{m,j}$.  The space
\[
        S_h:=\spanop\{B_1,\ldots,B_{n_B}\}
\]
has dimension $n_B$, consists of periodic splines of degree $m-1$ and global smoothness $C^{m-2}$, and satisfies $\sum_jB_j=1$.  Each $B_j$ has support of parameter length $mh$ and intersects at most $m+1$ intervals of the collocation mesh.  This parameter-based definition avoids any ordered-knot recursion with complex contour points across the parameter cut.

At the mesh points the interpolation matrix is circulant.  Its nonzero entries are
\begin{equation}\label{eq:collocationB}
\begin{array}{c|ccc}
 m & B_j(t_j^B) & B_j(t_{j-1}^B)=B_j(t_{j+1}^B) & \text{all other entries}\\
 \hline
 2 & 1 & 0 & 0\\
 3 & 3/4 & 1/8 & 0\\
 4 & 2/3 & 1/6 & 0.
\end{array}
\end{equation}
Thus the matrix is the identity for $m=2$ and is cyclic tridiagonal and strictly row diagonally dominant for $m=3,4$.  Its inverse is bounded independently of $n_B$.  Let $P_hv\in S_h$ denote the corresponding periodic spline interpolant defined by
\begin{equation}\label{eq:interpolation}
        (P_hv)(t_i^B)=v(t_i^B),
        \qquad i=1,\ldots,n_B.
\end{equation}

\begin{lemma}[Spline stability and approximation]\label{lem:splineapprox}
Let $0<\beta<\alpha\le1$.  There is a constant independent of $h$ such that
\begin{align}
        \norm{P_hv}_{\infty}
        &\le C\norm{v}_{\infty},\label{eq:Phsupstable}\\
        \norm{v-P_hv}_{\infty}
        &\le Ch^\alpha\norm{v}_{H_\alpha(\Ga)},\label{eq:Phsupapprox}\\
        \norm{P_hv}_{H_\beta(\Ga)}
        &\le C\norm{v}_{H_\beta(\Ga)},\label{eq:Phstable}\\
        \norm{v-P_hv}_{H_\beta(\Ga)}
        &\le Ch^{\alpha-\beta}\norm{v}_{H_\alpha(\Ga)}.
        \label{eq:Phapprox}
\end{align}
The last estimate also holds with $\alpha$ replaced by any
$\sigma\in(\beta,1]$.
\end{lemma}

\begin{proof}
These are standard stability and direct estimates for the interpolation rule \eqref{eq:interpolation} on a uniform periodic mesh; see \cite{deBoor,Schumaker}.  We give the ingredients needed below.  Write
\[
        P_hv=\sum_{j=1}^{n_B}a_jB_j,
        \qquad
        \mathbb B\bm a=\bm v,
        \qquad
        \bm v=(v(t_1^B),\ldots,v(t_{n_B}^B))^{\mathsf T}.
\]
For $m=2,3,4$, the circulant band matrix $\mathbb B$ is uniformly invertible.  In the tridiagonal cases, the entries of $\mathbb B^{-1}$ decay exponentially with the cyclic distance; hence
\[
        \norm{\bm a}_{\ell^\infty}\le C\norm{\bm v}_{\ell^\infty}
\]
and, for every cyclic shift $r$,
\[
        \max_j|a_{j+r}-a_j|
        \le C\max_j|v(t_{j+r}^B)-v(t_j^B)|.
\]
The partition of unity and local support of the periodic B-spline basis therefore imply \eqref{eq:Phsupstable} and
\[
        [P_hv]_{\beta,\Ga}\le C\norm{v}_{H_\beta(\Ga)},
\]
which proves \eqref{eq:Phstable} explicitly.

For $v\in H_\alpha(\Ga)$, the same coefficient bounds and polynomial reproduction give the standard local estimates
\[
        \norm{v-P_hv}_\infty\le Ch^\alpha\norm{v}_{H_\alpha},
        \qquad
        \norm{(P_hv)' }_\infty\le Ch^{\alpha-1}\norm{v}_{H_\alpha}.
\]
If $|t-s|\ge h$, the first estimate, divided by $|t-s|^\beta$, is bounded by $Ch^{\alpha-\beta}\norm{v}_{H_\alpha}$.  If $|t-s|<h$, the $\alpha$-H\"older estimate for $v$ and the derivative estimate for $P_hv$ give the same bound.  This proves \eqref{eq:Phsupapprox} and \eqref{eq:Phapprox}.  Replacing $\alpha$ by any $\sigma\in(\beta,1]$ gives the last assertion.
\end{proof}

\begin{remark}[Choice of the spline order]\label{rem:splineorder}
The method and the convergence proof are unchanged for $m=2,3,4$; only the basis values, support size, and matrix entries change.  Under the baseline assumption $u\in H_\alpha$ with $\alpha<1$, all three orders have the same guaranteed rate $h^{\alpha-\beta}$ because the data regularity, not the polynomial degree, is limiting.  Quadratic and cubic splines become advantageous when the continuous component has additional smoothness, whereas linear splines give the simplest and least expensive implementation.
\end{remark}

The enriched periodic spline space is
\begin{equation}\label{eq:enrichedspace}
        V_h^H
        :=S_h\oplus\spanop\{G_1,\ldots,G_{n_d}\}.
\end{equation}
The sum in \eqref{eq:enrichedspace} is direct: if
\(
 s_h+\sum_{j=1}^{n_d}c_jG_j=0
\)
with $s_h\in S_h$, then taking jumps at $t_i^d$ and using
\eqref{eq:jumpmatrix} gives $c_i=0$ for every $i$, and hence
$s_h=0$.

\begin{proposition}[Jackson estimate for the enriched space]\label{prop:jackson-enriched}
Let $0<\beta<\alpha<1$.  For every
$v\in X_\alpha=PH_\alpha(\Ga,\D)$, there exists
$v_h\in V_h^H$ having exactly the same jumps as $v$ such that
\begin{equation}\label{eq:jackson-enriched}
        \norm{v-v_h}_{X_\beta}
        \le C h^{\alpha-\beta}\norm{v}_{X_\alpha},
\end{equation}
where $C$ is independent of $h$ and $v$.  Consequently,
\begin{equation}\label{eq:best-jackson-enriched}
        \inf_{w_h\in V_h^H}
        \norm{v-w_h}_{X_\beta}
        \le C h^{\alpha-\beta}\norm{v}_{X_\alpha}.
\end{equation}
\end{proposition}

\begin{proof}
By Proposition~\ref{prop:decomposition}, write
\[
        v=v_C+J_v,
        \qquad
        J_v=\sum_{j=1}^{n_d}[v]_{t_j^d}G_j,
        \qquad
        v_C\in H_\alpha(\Ga),
\]
with
\(
 \norm{v_C}_{H_\alpha(\Ga)}\le C\norm{v}_{X_\alpha}.
\)
Choose
$
        v_h:=P_hv_C+J_v\in V_h^H.
$
Since $P_hv_C$ is globally continuous and $J_v$ reproduces the jump
vector of $v$, one has
\(
 [v_h]_{t_j^d}=[v]_{t_j^d}
\)
for every $j$.  Moreover, the jump components cancel exactly, so
$
        v-v_h=v_C-P_hv_C.
$
Using the continuous embedding of $H_\beta(\Ga)$ into $X_\beta$,
Lemma~\ref{lem:splineapprox}, and
\eqref{eq:decompositionbound}, we obtain
\[
        \norm{v-v_h}_{X_\beta}
        \le C\norm{v_C-P_hv_C}_{H_\beta(\Ga)}\\
        \le Ch^{\alpha-\beta}\norm{v_C}_{H_\alpha(\Ga)}
        \le Ch^{\alpha-\beta}\norm{v}_{X_\alpha}.
\]
This proves the enriched estimate \eqref{eq:jackson-enriched}; taking the infimum over
$V_h^H$ gives \eqref{eq:best-jackson-enriched}.
\end{proof}

The exact reduced collocation solution $u_h\in S_h$ is defined by
\begin{equation}\label{eq:pointcollocation}
        (Au_h)(t_i^B)=g(t_i^B),
        \qquad i=1,\ldots,n_B.
\end{equation}
Equivalently, the projected equation is
\begin{equation}\label{eq:projectedcollocation}
        A_hu_h=P_hg,
        \qquad
        A_h:=I-\lambda P_hK.
\end{equation}
The nodal collocation system \eqref{eq:pointcollocation} is equivalent to
\eqref{eq:projectedcollocation}.  Writing
$
        u_h(t)=\sum_{j=1}^{n_B}a_jB_j(t),
$
we obtain the dense linear system
\begin{equation}\label{eq:matrix}
        \sum_{j=1}^{n_B}\mathbb A_{ij}a_j=b_i,
        \qquad i=1,\ldots,n_B,
\end{equation}
where
\begin{align}
        \mathbb A_{ij}
        &=B_j(t_i^B)
          -\lambda\int_\Ga k(t_i^B,\tau)B_j(\tau)\,d\tau,
          \label{eq:matrixentries}\\
        b_i
        &=f_C(t_i^B)
          +\lambda\int_\Ga k(t_i^B,\tau)J_f(\tau)\,d\tau.
          \label{eq:rhsentries}
\end{align}
With the right-hand side entries \eqref{eq:rhsentries}, the final enriched approximation is
\begin{equation}\label{eq:finalapprox}
        \varphi_h=u_h+J_f.
\end{equation}

\section{Stability and convergence of exact collocation}\label{sec:convergence}

The smoothing property of $K$ makes the collocation stability a perturbation consequence of continuous invertibility rather than an independent unverified assumption.

\begin{lemma}[Operator-norm consistency]\label{lem:operatorconsistency}
Under Assumption~\ref{ass:kernel},
\begin{equation}\label{eq:operatorconsistency}
        \norm{(I-P_h)K}_{X_\beta\to X_\beta}
        \le Ch^{\mu-\beta}.
\end{equation}
\end{lemma}

\begin{proof}
For $v\in X_\beta$, Lemma~\ref{lem:Kmapping} gives
\[
        \norm{Kv}_{H_\mu}\le C\norm{v}_{X_\beta}.
\]
Apply Lemma~\ref{lem:splineapprox} with $\sigma=\mu$.
\end{proof}

\begin{theorem}[Uniform discrete stability]\label{thm:discretestability}
Under Assumptions~\ref{ass:kernel} and~\ref{ass:injectivity}, there exists $h_0=h_0(\lambda)>0$ such that, for $0<h\le h_0$, $A_h$ is invertible on $X_\beta$ and
\begin{equation}\label{eq:uniforminverse}
        \norm{A_h^{-1}}_{X_\beta\to X_\beta}\le C_\lambda.
\end{equation}
The constants $h_0$ and $C_\lambda$ may deteriorate as $\lambda$ approaches a characteristic value of $K$.
For every right-hand side in $S_h$, the solution of the $X_\beta$ equation belongs to $S_h$; hence \eqref{eq:projectedcollocation} and the nodal collocation system \eqref{eq:pointcollocation} have a unique solution.
\end{theorem}

\begin{proof}
By the operator consistency estimate \eqref{eq:operatorconsistency},
\[
        A_h-A=\lambda(I-P_h)K,
        \qquad
        \norm{A_h-A}_{X_\beta\to X_\beta}\to0.
\]
Since $A$ is invertible, choose $h_0$ so that
$
        \norm{A^{-1}(A_h-A)}\le\frac12.
$
The Banach--Neumann lemma yields invertibility and \eqref{eq:uniforminverse}.  If $A_hv=s_h\in S_h$, then
$
        v=s_h+\lambda P_hKv\in S_h.
$
Therefore the solution lies in the spline space.
\end{proof}

\begin{theorem}[Exact-collocation convergence]\label{thm:exactconvergence}
Let $0<\beta<\alpha<\mu\le1$, and suppose Assumptions~\ref{ass:kernel} and~\ref{ass:injectivity} hold.  For all sufficiently small $h$, the approximation \eqref{eq:finalapprox} exists uniquely and satisfies
\begin{equation}\label{eq:mainerror}
        \norm{\varphi-\varphi_h}_{X_\beta}
        \le C_\lambda h^{\alpha-\beta}\norm{f}_{X_\alpha}.
\end{equation}
The constant is independent of $h$ and $f$, but may depend on the fixed nonresonant parameter $\lambda$.
\end{theorem}

\begin{proof}
The jump components of $\varphi$ and $\varphi_h$ are identical, so $\varphi-\varphi_h=u-u_h$ is globally continuous.  The global and piecewise H\"older norms are equivalent on this subspace, and in particular
\[
        \norm{\varphi-\varphi_h}_{X_\beta}
        \le C\norm{u-u_h}_{H_\beta}.
\]
Because $Au=g$, one has
\[
        A_hu-P_hg=u-\lambda P_hKu-P_h(u-\lambda Ku)=(I-P_h)u.
\]
Together with $A_hu_h=P_hg$, this gives
\[
        A_h(u-u_h)=(I-P_h)u.
\]
Theorem~\ref{thm:discretestability}, Lemma~\ref{lem:splineapprox}, and estimate \eqref{eq:ubound} imply
\[
        \norm{u-u_h}_{H_\beta}
        \le C\norm{(I-P_h)u}_{H_\beta}
        \le Ch^{\alpha-\beta}\norm{u}_{H_\alpha}
        \le Ch^{\alpha-\beta}\norm{f}_{X_\alpha}.
\]
\end{proof}

\begin{remark}[Why a weaker H\"older exponent is used]\label{rem:nonsame}
Estimate \eqref{eq:mainerror} is stated in $X_\beta$ with $\beta<\alpha$.
Convergence in the same full H\"older exponent for every datum in
$X_\alpha$ generally requires additional smoothness.  For example,
$u\in H_{\alpha+\delta}$ yields a same-target estimate with the available
smoothness gap.
\end{remark}

\section{B-spline--Heaviside iterated extension}\label{sec:iterated}

The direct approximation
\[
        \varphi_h=u_h+J_f
\]
has exactly the prescribed physical jumps.  As in classical iterated
collocation \cite{AtkinsonGrahamSloan,Joe1985,GrahamJoeSloan1985}, it may be
inserted once into the original equation.

\begin{definition}[Iterated B-spline--Heaviside approximation]
\label{def:bsh-iterated}
For the exact-operator collocation solution $\varphi_h$, define
\begin{equation}\label{eq:bsh-iterated}
        \varphi_h^{\mathrm{it}}(t)
        :=f(t)+\lambda(K\varphi_h)(t),
        \qquad t\in\Ga.
\end{equation}
\end{definition}

Because $K\varphi_h$ is globally continuous,
\begin{equation}\label{eq:bsh-iterated-jumps}
        [\varphi_h^{\mathrm{it}}]_{t_j^d}
        =[f]_{t_j^d}
        =[\varphi]_{t_j^d},
        \qquad 1\le j\le n_d,
\end{equation}
so the iteration preserves the physical jumps and creates no mesh-interface
discontinuities.  Subtracting \eqref{eq:bsh-iterated} from the exact equation
gives
\begin{equation}\label{eq:bsh-iterated-error-identity}
        \varphi-\varphi_h^{\mathrm{it}}
        =\lambda K(\varphi-\varphi_h).
\end{equation}

\begin{theorem}[Improved raw and iterated estimates]
\label{thm:bsh-iterated}
Let
\[
        0<\beta<\alpha<\mu\le1,
\]
and suppose Assumptions~\ref{ass:kernel} and~\ref{ass:injectivity} hold.
Then, for all sufficiently small $h$, the exact B-spline--Heaviside
collocation solution satisfies
\begin{equation}\label{eq:bsh-raw-sup}
        \norm{\varphi-\varphi_h}_{\infty}
        \le
        C_\lambda h^\alpha\norm{f}_{X_\alpha},
\end{equation}
and its iterated extension satisfies
\begin{equation}\label{eq:bsh-iterated-bound}
        \norm{\varphi-\varphi_h^{\mathrm{it}}}_{H_\mu(\Ga)}
        \le
        C_\lambda h^\alpha\norm{f}_{X_\alpha}.
\end{equation}
Consequently, \eqref{eq:bsh-iterated-bound} also holds in $C(\Ga)$ and in every $H_\sigma(\Ga)$ with $0<\sigma\le\mu$.
\end{theorem}

\begin{proof}
First consider $A=I-\lambda K$ on $C(\Ga)$ with the uniform norm.
The map $K:C(\Ga)\to H_\mu(\Ga)$ is bounded and compact as a map into
$C(\Ga)$.  If $Av=0$ for some $v\in C(\Ga)$, then
$v=\lambda Kv\in H_\mu(\Ga)\subset X_\beta$, and Assumption
\ref{ass:injectivity} gives $v=0$.  Hence $A$ is boundedly invertible on
$C(\Ga)$.

By \eqref{eq:Phsupstable}, \eqref{eq:Phsupapprox}, and the smoothing property
of $K$,
\[
        \norm{(I-P_h)K}_{C(\Ga)\to C(\Ga)}
        \le Ch^\mu.
\]
Therefore $A_h=I-\lambda P_hK$ is uniformly invertible on $C(\Ga)$ for all
sufficiently small $h$.  Since
\[
        A_h(u-u_h)=(I-P_h)u,
\]
we obtain from \eqref{eq:Phsupapprox} and \eqref{eq:ubound}
\[
        \norm{u-u_h}_{\infty}
        \le
        C_\lambda h^\alpha\norm{u}_{H_\alpha}
        \le
        C_\lambda h^\alpha\norm{f}_{X_\alpha}.
\]
The exact and discrete jump components coincide, so
$\varphi-\varphi_h=u-u_h$, which proves \eqref{eq:bsh-raw-sup}.

Finally, \eqref{eq:bsh-iterated-error-identity} and the bound
$K:C(\Ga)\to H_\mu(\Ga)$ give
\[
        \norm{\varphi-\varphi_h^{\mathrm{it}}}_{H_\mu}
        \le
        |\lambda|\,C\norm{\varphi-\varphi_h}_{\infty}.
\]
Combining this estimate with \eqref{eq:bsh-raw-sup} proves
\eqref{eq:bsh-iterated-bound}.
\end{proof}

Theorem~\ref{thm:bsh-iterated} gives a strictly stronger rate for the
iterated approximation than the baseline $X_\beta$ estimate
$O(h^{\alpha-\beta})$.  It uses only the assumptions already imposed for
exact collocation and does not claim the higher superconvergence orders that
may occur for smoother data and special collocation points.

Although neither $\varphi$ nor $\varphi_h^{\mathrm{it}}$ is generally globally continuous when the prescribed jumps are nonzero, their difference belongs to $H_\mu(\Ga)$: the two functions have identical jumps, and \eqref{eq:bsh-iterated-error-identity} expresses the error as the image of the bounded raw error under the smoothing operator $K$.

For the fully discrete solution $\varphi_h^q$, the implementable extension is
\begin{equation}\label{eq:bsh-iterated-discrete}
        \varphi_h^{\mathrm{it},q}(t)
        :=f(t)+\lambda
        (K_{\mathrm{ext},h}^q\varphi_h^q)(t),
\end{equation}
where the source panels are split at spline breakpoints, physical jumps, and
any additional coefficient or cutoff breakpoints.  The quadrature conditions
needed to retain the exact-operator rates are stated in
Section~\ref{sec:quadrature}.

For
\[
        c\varphi-\lambda K\varphi=f,
        \qquad w=c\varphi,
\]
the raw method approximates $w$, and the iterated extension is
\begin{equation}\label{eq:bsh-iterated-variable}
        w_h^{\mathrm{it}}
        :=f+\lambda K\!\left(\frac{w_h}{c}\right),
        \qquad
        \varphi_h^{\mathrm{it}}
        :=\frac{w_h^{\mathrm{it}}}{c}.
\end{equation}

\section{Panelwise quadrature and the fully discrete method}\label{sec:quadrature}

In the periodic parameter,
\begin{equation}\label{eq:paramintegral}
        (Kv)(t)
        =
        \int_0^{2\pi}
        a(t,\theta)v(\gamma(\theta))\,d\theta,
        \qquad
        a(t,\theta):=
        k(t,\gamma(\theta))\gamma'(\theta).
\end{equation}
Let $\Pi_h$ be the common refinement of the periodic spline breakpoints, the
physical jump parameters, the cut, and any additional points at which a
coefficient or cutoff changes formula.  We assume
\begin{equation}\label{eq:panel-size}
        \max_{I\in\Pi_h}|I|\le C_\Pi h.
\end{equation}
On every panel, a spline $v_h\in S_h$ is a polynomial of degree at most
$m-1$, while the corrected jump component $J_f$ is affine in the parameter.

The following H\"older regularity assumption makes the quadrature error explicit.  It is required only for the algebraic fully discrete estimate;
the exact-collocation theory continues to use Assumption~\ref{ass:kernel}.

\begin{assumption}[Regularity in the integration variable]
\label{ass:quadrature-power}
There exist $0<\rho\le1$ and $L_\rho>0$ such that
\begin{equation}\label{eq:a-rho}
        \norm{a(\cdot,\theta)-a(\cdot,\eta)}_{H_\mu(\Ga)}
        \le
        L_\rho|\theta-\eta|^\rho,
        \qquad
        \theta,\eta\in[0,2\pi].
\end{equation}
\end{assumption}

For periodic analytic integrands, a composite trapezoidal rule may also be
attractive \cite{TrefethenWeideman}; however, the panelwise Gauss formulation
below remains applicable to the piecewise H\"older setting considered here.

Let $Q_I^q$ be the $q$-point Gauss--Legendre rule on a panel $I$, and define
\begin{equation}\label{eq:Khq}
        (K_h^qv)(t)
        :=
        \sum_{I\in\Pi_h}
        Q_I^q\!\left[
        a(t,\cdot)\,v(\gamma(\cdot))
        \right].
\end{equation}

\begin{proposition}[Explicit Gauss-panel consistency]
\label{prop:explicit-quadrature}
Suppose Assumptions~\ref{ass:kernel} and
\ref{ass:quadrature-power} hold, \eqref{eq:panel-size} is satisfied, and
\begin{equation}\label{eq:q-min}
        q\ge\left\lceil\frac m2\right\rceil.
\end{equation}
Then
\begin{align}
 \norm{(K-K_h^q)v_h}_{H_\mu(\Ga)}
 &\le
 Ch^\rho\norm{v_h}_{\infty},
 &&v_h\in S_h,\label{eq:q-explicit-spline}\\
 \norm{(K-K_h^q)J_f}_{H_\mu(\Ga)}
 &\le
 Ch^\rho\norm{J_f}_{\infty}.\label{eq:q-explicit-jump}
\end{align}
Consequently,
\begin{align}
 \norm{P_h(K-K_h^q)v_h}_{H_\beta(\Ga)}
 &\le
 Ch^\rho\norm{v_h}_{H_\beta(\Ga)},\label{eq:qop-explicit}\\
 \norm{P_h(K-K_h^q)J_f}_{H_\beta(\Ga)}
 &\le
 Ch^\rho\norm{J_f}_{\infty}.\label{eq:qrhs-explicit}
\end{align}
\end{proposition}

\begin{proof}
Fix a panel $I\in\Pi_h$ and choose $\theta_I\in I$.  Since
$v_h\circ\gamma$ is a polynomial of degree at most $m-1$ on $I$ and
$2q-1\ge m-1$, the Gauss rule integrates
$a(t,\theta_I)v_h(\gamma(\theta))$ exactly.  Hence the panel defect equals
the defect obtained after replacing $a(t,\theta)$ by
$a(t,\theta)-a(t,\theta_I)$.  Positivity of the Gauss weights and
\eqref{eq:a-rho} give
\[
 \left\|
 \int_I av_h-Q_I^q(av_h)
 \right\|_{H_\mu}
 \le
 C|I|^{1+\rho}\norm{v_h}_{\infty,I}.
\]
Summing over the panels and using
\[
        \sum_{I\in\Pi_h}|I|^{1+\rho}
        \le
        \left(\max_{I\in\Pi_h}|I|\right)^\rho
        \sum_{I\in\Pi_h}|I|
        \le Ch^\rho
\]
proves \eqref{eq:q-explicit-spline}.  On each jump-aligned panel,
$J_f\circ\gamma$ is affine; condition \eqref{eq:q-min} also integrates the
frozen-kernel affine term exactly, and the same argument proves
\eqref{eq:q-explicit-jump}.  The right-hand estimate follows analogously.  Finally, use the embedding $H_\mu\hookrightarrow H_\beta$ and the $H_\beta$ stability of $P_h$.
\end{proof}

The fully discrete spline component $u_h^q\in S_h$ is defined by
\begin{equation}\label{eq:fullydiscrete}
        \left(I-\lambda P_hK_h^q\right)u_h^q
        =
        P_h\left(f_C+\lambda K_h^qJ_f\right),
\end{equation}
and $\varphi_h^q=u_h^q+J_f$.

\begin{theorem}[Explicit fully discrete convergence]
\label{thm:quadrature}
Let
\[
        0<\beta<\alpha<\mu\le1.
\]
Under the hypotheses of Theorem~\ref{thm:exactconvergence},
Assumption~\ref{ass:quadrature-power}, \eqref{eq:panel-size}, and
\eqref{eq:q-min}, the fully discrete problem is uniquely solvable for all
sufficiently small $h$ and
\begin{equation}\label{eq:qerror}
        \norm{\varphi-\varphi_h^q}_{X_\beta}
        \le
        C_\lambda
        \left(h^{\alpha-\beta}+h^\rho\right)
        \norm{f}_{X_\alpha}.
\end{equation}
In particular, the exact-collocation rate is preserved whenever
\begin{equation}\label{eq:rho-condition}
        \rho\ge\alpha-\beta.
\end{equation}
\end{theorem}

\begin{proof}
Let
\[
        A_{h,S}:=\left.(I-\lambda P_hK)\right|_{S_h}:S_h\to S_h,
        \qquad
        A_{h,S}^q:=\left.(I-\lambda P_hK_h^q)\right|_{S_h}:S_h\to S_h.
\]
Equip $S_h$ with the norm inherited from $H_\beta(\Ga)$.  Theorem~\ref{thm:discretestability} implies that $A_{h,S}^{-1}$ is uniformly bounded, while Proposition~\ref{prop:explicit-quadrature} gives
\[
        \norm{A_{h,S}^q-A_{h,S}}_{S_h\to S_h}
        \le C|\lambda|h^\rho.
\]
The Banach--Neumann lemma therefore yields uniform invertibility of $A_{h,S}^q$ for sufficiently small $h$.  Subtracting the exact and fully discrete equations yields
\[
 \left(I-\lambda P_hK_h^q\right)(u_h-u_h^q)
 =
 \lambda P_h(K-K_h^q)u_h
 +
 \lambda P_h(K-K_h^q)J_f.
\]
The exact discrete solution is uniformly bounded, and the jump-component
estimate following from Proposition~\ref{prop:decomposition} gives
$\norm{J_f}_\infty\le C\norm{f}_{X_\alpha}$.  Therefore
\[
        \norm{u_h-u_h^q}_{H_\beta}
        \le
        C_\lambda h^\rho\norm{f}_{X_\alpha}.
\]
Combining this with Theorem~\ref{thm:exactconvergence} proves
\eqref{eq:qerror}.
\end{proof}

\begin{remark}[Choice of the Gauss order]\label{rem:q-choice}
No growth of $q$ with $h^{-1}$ is required under
Assumption~\ref{ass:quadrature-power}: a fixed
$q\ge\lceil m/2\rceil$ is sufficient, provided
\eqref{eq:rho-condition} holds.  For $m=2,3,4$, the minimum values are
$q=1,2,2$, respectively; the order $q=10$ used in the experiments is
therefore conservative.  If only a general modulus of continuity is known,
the same proof gives a defect of order $\omega_a(C_\Pi h)$ instead of
$h^\rho$.

In computations one may also increase $q$ until
\[
 D_h(q):=
 \max\left\{
 \frac{\norm{\mathbb A_h^q-\mathbb A_h^{q_{\rm ref}}}_\infty}
      {\norm{\mathbb A_h^{q_{\rm ref}}}_\infty},
 \frac{\norm{b_h^q-b_h^{q_{\rm ref}}}_\infty}
      {\norm{b_h^{q_{\rm ref}}}_\infty}
 \right\}
 \le
 \vartheta h^{\alpha-\beta},
 \qquad 0<\vartheta<1,
\]
with a higher reference order $q_{\rm ref}$.  This is a practical
verification of the theoretical scale, not an additional assumption in
Theorem~\ref{thm:quadrature}.
\end{remark}

\section{A piecewise H\"older leading coefficient}\label{sec:variablecoefficient}

Consider the more general equation
\begin{equation}\label{eq:variablec}
        c(t)\varphi(t)-\lambda(K\varphi)(t)=f(t),
        \qquad t\in\Ga,
\end{equation}
where $c\in X_\alpha$ has the same prescribed jump set as $f$.  Assume
\begin{equation}\label{eq:cnonzero}
        \inf_{t\in\Ga\setminus\D}|c(t)|\ge c_0>0
\end{equation}
and that all lateral values of $c$ are nonzero.  Under the nonvanishing condition \eqref{eq:cnonzero}, equation \eqref{eq:variablec} is a Fredholm equation of the second kind after multiplication by $c^{-1}$.  Direct division, however, produces a kernel coefficient $c(t)^{-1}$ that is discontinuous in the observation variable and obscures the jump structure.  A better transformation is
\begin{equation}\label{eq:wtransform}
        w:=c\varphi.
\end{equation}
Under the transformation \eqref{eq:wtransform}, equation \eqref{eq:variablec} becomes
\begin{equation}\label{eq:wequation}
        w-\lambda\widetilde K w=f,
        \qquad
        (\widetilde K w)(t)
        :=\int_\Ga k(t,\tau)c(\tau)^{-1}w(\tau)\,d\tau.
\end{equation}
Because $c^{-1}\in X_\alpha$ and multiplication by $c^{-1}$ is bounded on $X_\beta$, the bounds \eqref{eq:kbound}--\eqref{eq:kholder} imply
\[
        \widetilde K:X_\beta\longrightarrow H_\mu(\Ga)
\]
is bounded and compact as an operator on $X_\beta$.  The transformed kernel $k(t,\tau)c(\tau)^{-1}$ is generally only piecewise uniformly continuous in $\tau$, because $c^{-1}$ may jump at $\D$.  This does not affect the mapping, compactness, jump-transfer, or exact-collocation arguments, which use only boundedness and the uniform H\"older estimate in $t$.  For the fully discrete method, the quadrature panels must additionally be split at the points of $\D$.  If Assumption~\ref{ass:quadrature-power} holds with exponent $\rho$, then on each continuity arc $c^{-1}$ is $\alpha$-H\"older and the transformed parametrized density satisfies the same condition with exponent $\widetilde\rho:=\min\{\rho,\alpha\}$.  Hence Theorem~\ref{thm:quadrature} applies to the transformed equation with a quadrature term $O(h^{\widetilde\rho})$.  Since $\alpha>\alpha-\beta$, the condition $\rho\ge\alpha-\beta$ remains sufficient to preserve the exact-collocation rate.

\begin{assumption}[Transformed nonresonance]\label{ass:transformed-nonresonance}
For the operator $\widetilde K$ in \eqref{eq:wequation},
\[
        \ker(I-\lambda\widetilde K)=\{0\}
        \qquad\text{in }X_\beta.
\]
Equivalently, if $\lambda\ne0$, then
$\lambda^{-1}\notin\sigma(\widetilde K:X_\beta\to X_\beta)$.
\end{assumption}

A sufficient, but not necessary, condition for Assumption~\ref{ass:transformed-nonresonance} is $|\lambda|\,\norm{\widetilde K}_{X_\beta\to X_\beta}<1$.  Under Assumption~\ref{ass:transformed-nonresonance}, the operator-theoretic and exact-collocation results of the preceding sections apply to $w$, while the panelwise quadrature result applies with the additional splitting just described.  In particular, the jump identity becomes
\begin{equation}\label{eq:wjumps}
        [w]_{t_j^d}=[f]_{t_j^d},
\end{equation}
and the identity \eqref{eq:wjumps} retains the triangular jump reconstruction for $w$, even though it is generally false for $\varphi$ itself.

Indeed, if $q_j:=(K\varphi)(t_j^d)$, continuity of $K\varphi$ gives
\begin{equation}\label{eq:philateralvariablec}
 \varphi(t_j^d\pm0)
 =\frac{f(t_j^d\pm0)+\lambda q_j}{c(t_j^d\pm0)},
\end{equation}
and therefore
\begin{equation}\label{eq:phijumpvariablec}
 [\varphi]_{t_j^d}
 =\frac{f(t_j^d+0)}{c(t_j^d+0)}
  -\frac{f(t_j^d-0)}{c(t_j^d-0)}
  +\lambda q_j
  \left(\frac1{c(t_j^d+0)}-\frac1{c(t_j^d-0)}\right).
\end{equation}
Relations \eqref{eq:philateralvariablec} and \eqref{eq:phijumpvariablec} show that the jump of $\varphi$ is not determined by the data alone when $c$ is discontinuous.

The practical method is to solve \eqref{eq:wequation} by the previous B-spline--Heaviside construction and then use the recovery formula
\begin{equation}\label{eq:phifromw}
        \varphi_h=c^{-1}w_h.
\end{equation}
Formula \eqref{eq:phifromw} reconstructs the approximation of the original unknown.  The resulting matrix and right-hand side are
\begin{align}
 \mathbb A_{ij}^{(c)}
 &=B_j(t_i^B)
 -\lambda\int_\Ga k(t_i^B,\tau)c(\tau)^{-1}B_j(\tau)\,d\tau,
 \label{eq:matrixc}\\
 b_i^{(c)}
 &=f_C(t_i^B)
 +\lambda\int_\Ga k(t_i^B,\tau)c(\tau)^{-1}J_f(\tau)\,d\tau.
 \label{eq:rhsc}
\end{align}
The entries \eqref{eq:matrixc}--\eqref{eq:rhsc} differ from the constant-coefficient case only by the factor $c(\tau)^{-1}$.  Since multiplication by $c^{-1}$ is bounded on $X_\beta$, the error estimate for $w_h$ immediately yields
\begin{equation}\label{eq:variablecerror}
        \norm{\varphi-\varphi_h}_{X_\beta}
        \le Ch^{\alpha-\beta}\norm{f}_{X_\alpha}.
\end{equation}
Estimate \eqref{eq:variablecerror} assumes exact knowledge of $c$ and its inverse.  If $c$ vanishes or has a zero lateral value, this reduction is unavailable and the problem is no longer covered by the present second-kind theory.

\section{Implementation algorithm and efficiency}\label{sec:algorithm}

\begin{algorithm}[Periodic B-spline--Heaviside collocation]\label{alg:direct}
Given $\gamma$, $k$, $\lambda$, the jump set $\D$, and either exact lateral data or samples of $f$, choose $m\in\{2,3,4\}$ and perform the following steps.
\begin{enumerate}[label=\textbf{Step \arabic*.},leftmargin=2.6cm]
\item Choose a cut $t_*\notin\D$, shift the periodic parametrization so that $\gamma(0)=t_*$, and order the jump parameters $\theta_j^d$.

\item If the lateral values are given, compute
\[
        \delta_j=f(t_j^d+0)-f(t_j^d-0).
\]
If only samples are available, apply Algorithm~\ref{alg:lateral} and use $\widehat\delta_j$.  Construct $G_j=H_j-R$, the jump component $J_f$, and the continuous data component $f_C=f-J_f$ (or their reconstructed counterparts).

\item Select $n_B$, set $h=2\pi/n_B$, and choose a phase satisfying \eqref{eq:phase}.  Construct the periodic order-$m$ basis \eqref{eq:periodicB}.

\item Form the common panel refinement from the spline breakpoints, the jump parameters, and the cut.  Under Assumption~\ref{ass:quadrature-power}, choose a fixed Gauss order
\[
        q\ge\left\lceil\frac m2\right\rceil.
\]
For each target $t_i^B$, compute
\[
        r_i^J=(K_h^qJ_f)(t_i^B)
\]
panel by panel.

\item Set $b_i$ according to \eqref{eq:rhsentries}, that is,
\[
        b_i=f_C(t_i^B)+\lambda r_i^J.
\]

\item Assemble \eqref{eq:matrixentries}.  Because $B_j$ has compact support of length $mh$, each integral column requires quadrature only on the intersected local panels; no global integration of an individual spline basis function is needed.

\item Solve the dense system \eqref{eq:matrix} and define
\[
        u_h=\sum_{j=1}^{n_B}a_jB_j,
        \qquad \varphi_h=u_h+J_f.
\]

\item When an iterated approximation is required, evaluate
\[
        \varphi_h^{\mathrm{it},q}(t)
        =f(t)+\lambda(K_{\mathrm{ext},h}^q\varphi_h^q)(t)
\]
at the requested target points using the panelwise source partition.  For the
variable-coefficient problem, apply \eqref{eq:bsh-iterated-variable} to the
transformed unknown $w_h$.

\item Verify the nodal residual and the recovered jumps
\[
        [\varphi_h]_{t_j^d}
        =[\varphi_h^{\mathrm{it},q}]_{t_j^d}
        =\delta_j.
\]
When an a posteriori quadrature check is desired, increase $q$ until the
defect $D_h(q)$ from Remark~\ref{rem:q-choice} is below the target scale
$\vartheta h^{\alpha-\beta}$.
\end{enumerate}
\end{algorithm}

\subsection*{Operation counts}
For a $q$-point rule on each active cell and fixed $m\in\{2,3,4\}$, the spline part of the matrix requires $O((m+1)qn_B^2)=O(qn_B^2)$ kernel evaluations because each of the $n_B$ basis functions intersects at most $m+1$ collocation intervals and is evaluated at $n_B$ targets.  Under the panel-size condition \eqref{eq:panel-size} required by Proposition~\ref{prop:explicit-quadrature}, the jump contribution is also integrated on $N_p=O(n_B)$ panels and therefore requires $O(qN_pn_B)=O(qn_B^2)$ kernel evaluations at all collocation targets.  A lower cost for this contribution is possible only with a separate global, adaptive, or higher-order quadrature whose tolerance is chosen consistently with the target discretization error.  The matrix is dense because the Fredholm kernel is nonlocal; a direct factorization costs $O(n_B^3)$ operations and $O(n_B^2)$ memory.  If several right-hand sides share the same kernel, contour, and mesh, the matrix and its factorization can be reused, while only the jump reconstruction and right-hand side must be recomputed.

\section{Numerical experiments}\label{sec:numerics}

This section examines the phase-robust convergence rate in \eqref{eq:mainerror} predicted by
Theorem~\ref{thm:exactconvergence}, the need for the Heaviside enrichment, the
behavior of the three spline orders, the iterated extension of
Section~\ref{sec:iterated}, the variable-coefficient reduction of
Section~\ref{sec:variablecoefficient}, the one-sided reconstruction
procedure of Section~\ref{sec:lateral}, and the comparison with established
discontinuity-adapted alternatives.  All computations were performed in
MATLAB R2025b with complex arithmetic using the accompanying implementation
and test validation package~\cite{CapceleaCapceleaSoftware2026}.

\subsection{Experimental setting and robust phase selection}

We use the original contour parameter $\xi\in[0,2\pi)$ to specify the
geometry and then introduce the computational parameter
\[
        \theta=\xi-\xi_*\pmod{2\pi},
\]
where the auxiliary cut $\xi_*$ is chosen outside the physical jump set.
Thus $\theta=0$ always represents the cut point $t_*$, even in the tests in
which a physical jump corresponds to the original value $\xi=0$.

For Tests~1 and 2, the exact solution is manufactured in the form
\begin{equation}\label{eq:numerical-exact-solution}
 \varphi(\theta)=u_\alpha(\theta)
       +\sum_{j=1}^{n_d}\delta_jG_j(\theta),
\end{equation}
where
\begin{align}
 u_\alpha(\theta)
  &=(0.55+0.08i)+A_c d_{\mathbb T}(\theta,\theta_c)^\alpha
       +A_s q(\theta),\label{eq:numerical-continuous-part}\\
 d_{\mathbb T}(\theta,\eta)
  &:=2\left|\sin\frac{\theta-\eta}{2}\right|,\qquad
 q(\theta):=\cos\theta+0.45i\sin(2\theta)-0.20\cos(3\theta).
\end{align}
The right-hand side is evaluated from the exact state by
$f=\varphi-\lambda K\varphi$.  Test~3 instead manufactures
$w=c\varphi$ in the form \eqref{eq:numerical-exact-solution}, sets
$\varphi=w/c$, and evaluates $f=w-\lambda K\varphi$.  Test~4 retains the
jump decomposition \eqref{eq:numerical-exact-solution} but replaces
\eqref{eq:numerical-continuous-part} by the localized cusp described below.
Consequently, all solutions and right-hand sides are genuinely complex-valued.

The mesh levels are
$
        n_B=32,64,128,256,512,1024,
$
and all three orders $m=2,3,4$ are tested.  For each pair $(m,n_B)$, three
mesh phases $\rho$ are selected subject to the uniform separation condition
\begin{equation}\label{eq:numerical-phase-condition}
 \min_{1\le j\le n_d}
 \frac{\operatorname{dist}(\theta_j^d,\Theta_{n_B}(\rho))}{h}
 \ge 0.08,
 \qquad
 \Theta_{n_B}(\rho)=\{\rho+\ell h\pmod{2\pi}:0\le\ell<n_B\}.
\end{equation}
Condition \eqref{eq:numerical-phase-condition} prevents a collocation node from approaching a prescribed jump too closely.  The tables report the median over the three phases, while the maximum over
the same phases is used as a conservative robustness indicator.  Observed
rates are least-squares slopes of the median errors on the log--log scale for
$n_B\ge64$; hence they are less sensitive than two-level quotients to the
relative position of a H\"older cusp and the spline grid.

The collocation integrals are evaluated by panelwise Gauss--Legendre
quadrature of order $10$, with the panels split at spline knots, physical
jump points, H\"older cusp locations, and the additional coefficient or
cutoff breakpoints when they are present.  A separate order-$18$ rule is used for the quadrature-defect
check, whereas the manufactured right-hand side is generated with an
independent order-$120$ panel rule.  The reported errors are
\[
 E_\infty:=\|\varphi-\varphi_h\|_\infty,
 \qquad
 E_{PH_\beta}:=\|\varphi-\varphi_h\|_{PH_\beta,\mathrm{disc}}.
\]
The discrete $PH_\beta$ norm is evaluated on a grid of at least $4096$
points, separately on every continuity arc, so that no difference quotient
crosses a jump.

\subsection{Discontinuity-adapted comparison methods}
\label{sec:alternative-comparison}

The numerical comparison is separated into raw approximations and one-step
extensions.  The raw group consists of BSH, arcwise spline collocation in the
Atkinson--Graham--Sloan framework \cite{AtkinsonGrahamSloan}, and Joe's
discontinuous degree-$(m-1)$ Gauss collocation \cite{Joe1985}.  The iterated
group consists of the corresponding three extensions
$f+\lambda K\varphi_h$ and a panelwise Gauss--Nystr\"om method
\cite{Atkinson,JoeDiscrete1985}.  This separation avoids comparing a direct
trial-space approximation with a postprocessed one.

All alternatives are adapted in the periodic parameter by aligning their
panels with the physical jumps.  AGS uses independent open-knot spline
spaces on the continuity arcs and distinct lateral endpoint values.  Joe
uses local Lagrange polynomials at the $m$ Gauss points of each panel without
continuity constraints; consequently, its raw approximation may have
nonphysical panel-interface jumps.  Nystr\"om uses the same jump-aligned
Gauss panels and its standard extension formula.  For Test~3 the methods are
applied to $w-\lambda K(w/c)=f$ and then divided by $c$.

The comparison uses Tests~1--4, $m=2,3,4$, and
$n_B=32,\ldots,1024$, with closely matched solved dimensions and a common
independent validation grid.  The article reports only representative
final-level data and aggregate conclusions.  Complete all-order tables,
local rates, jump diagnostics, figures, and machine-readable results are
provided in the accompanying Zenodo archive
\cite{CapceleaCapceleaSoftware2026}.

\subsection{Tests 1 and 2: piecewise H\"older data on noncircular contours}

In Test~1,
\begin{align*}
 \gamma_1(\xi)
 &=1.34\cos\xi+0.79i\sin\xi+0.065e^{3i\xi},\\
 \xi_*&=0.36\pi,
 \qquad
 (\theta_1^d,\theta_2^d)=(0.52\pi,1.64\pi),\\
 (\delta_1,\delta_2)
 &=(-0.07+0.09i,\ 0.29+0.04i),
\end{align*}
while $\alpha=0.87$, $\theta_c=1.03\pi$,
$A_c=0.16+0.10i$, and $A_s=0.12+0.03i$.  We take
\[
 \lambda=0.66,
 \qquad
 k_1(t,\tau)=0.05\bigl(\cos(0.24(t-\tau))+0.06t\tau\bigr).
\]
The first physical jump occurs at the original parameter $\xi=0$, whereas
the cut is at $\xi_*=0.36\pi$, so the corrected generators are essential for
avoiding a cut artifact.

Test~2 uses the star-shaped polar contour
\[
\begin{gathered}
\begin{aligned}
r(\xi)&=1+0.11\cos(3\xi)+0.05\sin(2\xi),
\quad
\gamma_2(\xi)=r(\xi)e^{i\xi},\\
\xi_*&=1.76\pi,
\quad
(\theta_1^d,\theta_2^d,\theta_3^d)
=(0.24\pi,0.83\pi,1.58\pi),
\end{aligned}\\[2mm]
(\delta_1,\delta_2,\delta_3)
=(0.24-0.07i,\,-0.08+0.11i,\,0.06+0.02i).
\end{gathered}
\]
Here $\alpha=0.79$, $\theta_c=1.21\pi$,
$A_c=0.20-0.06i$, $A_s=0.09+0.07i$, and
\[
 \lambda=0.80,
 \qquad
 k_2(t,\tau)
 =0.04\frac{1+0.08t+0.05\tau}{3-0.12t\overline\tau}.
\]
The jumps are unbalanced, since $\sum_j\delta_j\ne0$, and the first one
again corresponds to the original parameter $\xi=0$ rather than to the cut.

Table~\ref{tab:test12-convergence} gives the phase-median errors for cubic
B-splines.  The parenthesized values are the phase maxima in the principal
piecewise H\"older norm.  The fitted $PH_\beta$ rates are $0.451$ in Test~1
and $0.465$ in Test~2, compared with the guaranteed values
$\alpha-\beta=0.32$ and $0.34$, respectively.  At $n_B=1024$, the ratios
between the classical continuous-spline error and the enriched error in the
principal $PH_\beta$ norm are approximately $279$ and $92$.  In these two tests, a small collocation residual for the continuous method does not prevent a large error in the piecewise H\"older norm, because the trial space cannot reproduce the prescribed jumps.

\begin{table}[t]
\centering
\caption{Phase-robust convergence for Tests~1 and 2 with $m=4$.  The
$PH_\beta$ entry is the median, followed by the phase maximum in parentheses.}
\label{tab:test12-convergence}
\small
\begin{tabular}{rccccc}
\toprule
&\multicolumn{2}{c}{Test 1: $\beta=0.55$}
&\multicolumn{2}{c}{Test 2: $\beta=0.45$}\\
\cmidrule(lr){2-3}\cmidrule(lr){4-5}
$n_B$ & $E_\infty$ & $E_{PH_\beta}$
      & $E_\infty$ & $E_{PH_\beta}$\\
\midrule
32   & $8.64\!\times\!10^{-3}$ & $6.77\!\times\!10^{-2}\ (8.10\!\times\!10^{-2})$
     & $1.54\!\times\!10^{-2}$ & $8.95\!\times\!10^{-2}\ (1.01\!\times\!10^{-1})$\\
64   & $6.41\!\times\!10^{-3}$ & $5.70\!\times\!10^{-2}\ (6.14\!\times\!10^{-2})$
     & $1.32\!\times\!10^{-2}$ & $7.56\!\times\!10^{-2}\ (7.61\!\times\!10^{-2})$\\
128  & $5.14\!\times\!10^{-3}$ & $4.58\!\times\!10^{-2}\ (4.66\!\times\!10^{-2})$
     & $5.35\!\times\!10^{-3}$ & $5.09\!\times\!10^{-2}\ (5.65\!\times\!10^{-2})$\\
256  & $2.60\!\times\!10^{-3}$ & $3.39\!\times\!10^{-2}\ (3.53\!\times\!10^{-2})$
     & $4.24\!\times\!10^{-3}$ & $4.07\!\times\!10^{-2}\ (4.09\!\times\!10^{-2})$\\
512  & $8.44\!\times\!10^{-4}$ & $2.32\!\times\!10^{-2}\ (2.54\!\times\!10^{-2})$
     & $1.35\!\times\!10^{-3}$ & $2.66\!\times\!10^{-2}\ (3.21\!\times\!10^{-2})$\\
1024 & $4.40\!\times\!10^{-4}$ & $1.68\!\times\!10^{-2}\ (1.92\!\times\!10^{-2})$
     & $6.89\!\times\!10^{-4}$ & $2.09\!\times\!10^{-2}\ (2.54\!\times\!10^{-2})$\\
\bottomrule
\end{tabular}
\end{table}

\begin{figure}[t]
\centering
\includegraphics[width=0.94\textwidth]{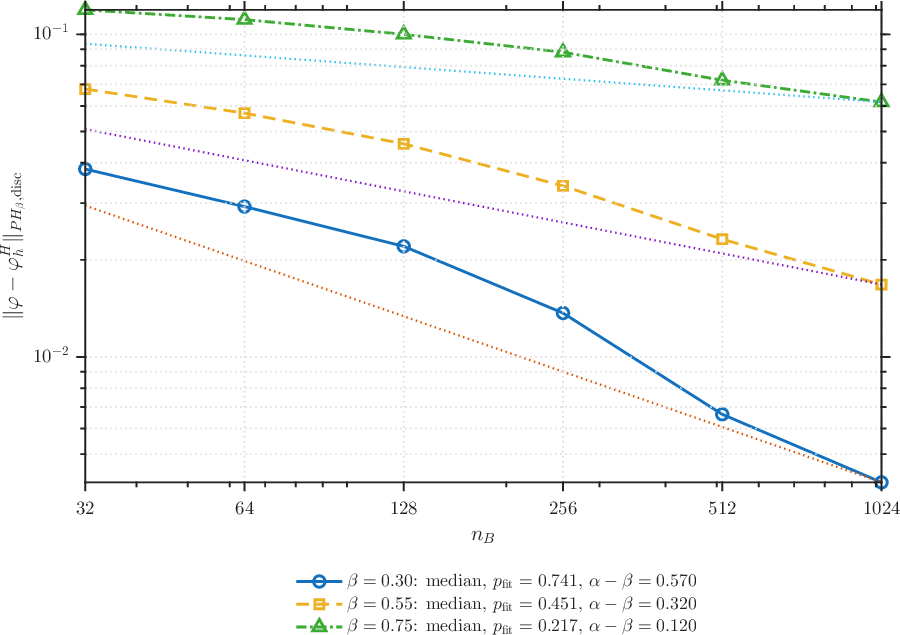}
\caption{Test~1, $m=4$: phase-median discrete $PH_\beta$ errors and reference
slopes $n_B^{-(\alpha-\beta)}$.}
\label{fig:test1-convergence}
\end{figure}

\begin{figure}[t]
\centering
\includegraphics[width=0.94\textwidth]{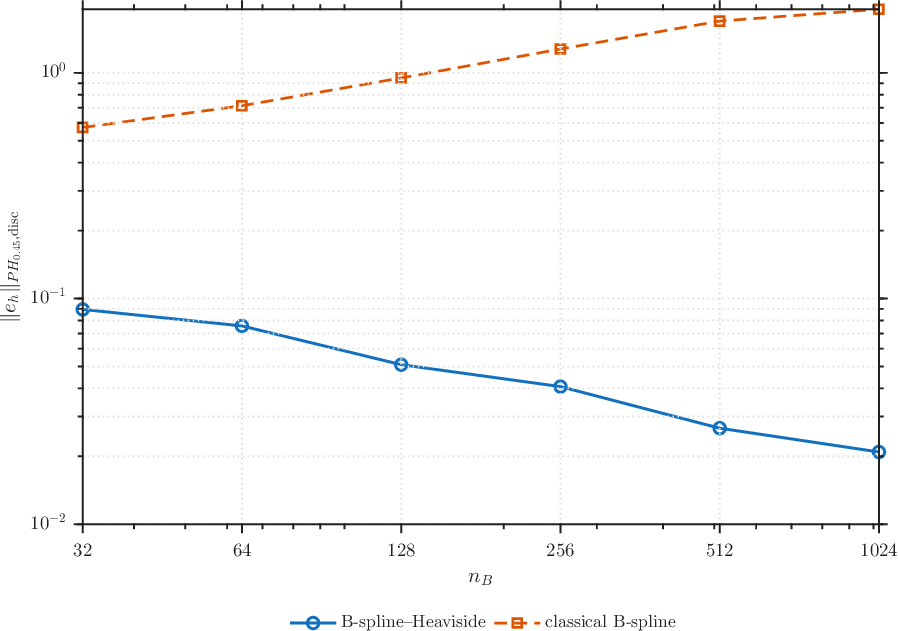}
\caption{Test~2, $m=4$: comparison of the phase-median $PH_{0.45}$ errors for
the B-spline--Heaviside method and the classical continuous B-spline method.}
\label{fig:test2-comparison}
\end{figure}

Figure~\ref{fig:test1-convergence} shows median errors for three values of $\beta$ that are compatible with the reference log--log slopes predicted by \eqref{eq:mainerror}.  The direct comparison
in Figure~\ref{fig:test2-comparison} is even more diagnostic: the enriched error decreases, whereas in this experiment the continuous-spline $PH_{0.45}$ error grows as the mesh represents the unresolved jump more sharply.

The continuous B-spline calculation is used here as a structural ablation baseline, not as a claim of superiority over every method designed for discontinuous solutions.  A separate comparison with arcwise AGS collocation, Joe's discontinuous Gauss collocation, their iterated extensions, and the Gauss--Nystr\"om method is specified in Subsection~\ref{sec:alternative-comparison}.

\subsection{Influence of the B-spline order \texorpdfstring{$m=2,3,4$}{m=2,3,4}}

Table~\ref{tab:order-comparison} compares the three orders in Test~2.  The
reported time is the median wall-clock time for the complete final level
$n_B=1024$ and is included only as a relative implementation indicator.
Linear splines give the smallest condition number, while quadratic and cubic
splines reduce the error constant.  The improvement from $m=3$ to $m=4$ is
small, which is expected because the global approximation rate is limited by
the H\"older regularity rather than by the polynomial degree.  The same trend
was observed in Tests~1 and 3.

\begin{table}[t]
\centering
\caption{Influence of the spline order in Test~2 at $n_B=1024$.  The fitted
rate is computed from the phase-median $PH_{0.45}$ errors for $n_B\ge64$.}
\label{tab:order-comparison}
\small
\begin{tabular}{cccccc}
\toprule
$m$ & $E_\infty$ & $E_{PH_{0.45}}$ & fitted rate
    & $\max\operatorname{cond}(\mathbb A_h)$ & time (s)\\
\midrule
2 & $8.53\times10^{-4}$ & $2.43\times10^{-2}$ & $0.489$ & $1.072$ & $0.777$\\
3 & $7.22\times10^{-4}$ & $2.12\times10^{-2}$ & $0.468$ & $2.071$ & $0.812$\\
4 & $6.89\times10^{-4}$ & $2.09\times10^{-2}$ & $0.465$ & $3.106$ & $0.862$\\
\bottomrule
\end{tabular}
\end{table}

For all three orders, the matrix condition numbers remain essentially
independent of $n_B$.  Across Tests~1--3, their maxima lie below $1.26$ for
$m=2$, $2.24$ for $m=3$, and $3.36$ for $m=4$.  Over the tested meshes, this behavior is consistent with the uniform stability mechanism of Section~\ref{sec:convergence}; finite computations alone do not constitute a proof of mesh-uniform stability.

\subsection{Test 3: a complex nonvanishing leading coefficient}

Test~3 examines the reduction of Section~\ref{sec:variablecoefficient}.  The
contour is
\[
 r(\xi)=1+0.14\cos(2\xi)-0.05\sin(3\xi)+0.03\cos\xi,
 \qquad \gamma_3(\xi)=r(\xi)e^{i\xi},
\]
with $\xi_*=0.14\pi$ and
\[
 (\theta_1^d,\theta_2^d,\theta_3^d)
 =(0.17\pi,0.89\pi,1.55\pi).
\]
The prescribed jumps of $w=c\varphi$ are
\[
 (0.23+0.05i,\ -0.12+0.08i,\ 0.07-0.03i),
\]
and the continuous component has $\alpha=0.81$,
$\theta_c=0.62\pi$, $A_c=0.19+0.07i$, and $A_s=0.10-0.04i$.
The leading coefficient is
\begin{equation}\label{eq:test3-c}
 c(\theta)=c_j^0+(0.11+0.04i)
 d_{\mathbb T}(\theta,1.13\pi)^{0.84}
 +(0.035-0.025i)\cos\theta.
\end{equation}
The coefficient in \eqref{eq:test3-c} is used for $\theta$ on the outgoing arc $(\theta_j^d,\theta_{j+1}^d)$, with cyclic
indexing and
\[
 (c_1^0,c_2^0,c_3^0)
 =(1.42+0.18i,\ 1.03-0.16i,\ 1.24+0.11i).
\]
Finally,
\[
 \lambda=0.70,
 \qquad
 k_3(t,\tau)=0.042e^{0.10t}
 \bigl(1+0.07\tau+0.04\overline\tau+0.03t\tau\bigr).
\]
Since $0\le d_{\mathbb T}\le2$, the reverse triangle inequality gives the analytic bound
\[
 \inf_\theta|c(\theta)|
 \ge
 \min_j|c_j^0|
 -|0.11+0.04i|\,2^{0.84}
 -|0.035-0.025i|
 >0.789.
\]
Thus the coefficient is rigorously separated from zero.  A dense diagnostic grid gives
\[
 \min_{\theta}|c(\theta)|=1.007,
 \qquad
 \max_{\theta}|\operatorname{Im}c(\theta)|=0.243,
\]
and confirms that the coefficient is genuinely complex.
The corresponding exact solution and right-hand side also have nonzero
imaginary parts, with maxima $0.409$ and $0.308$, respectively.

\begin{figure}[t]
\centering
\includegraphics[width=0.94\textwidth]{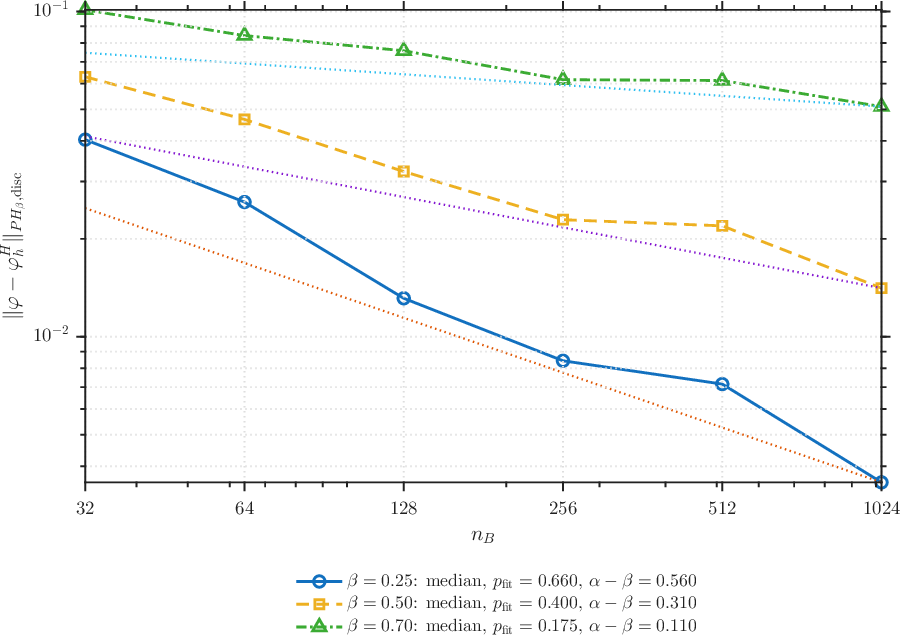}
\caption{Test~3, $m=4$: phase-median discrete $PH_\beta$ errors for the
transformed equation $w-\lambda\widetilde K w=f$.}
\label{fig:test3-convergence}
\end{figure}

For $m=4$ and $\beta=0.50$, the phase-median $PH_\beta$ error decreases from
$6.30\times10^{-2}$ at $n_B=32$ to $1.41\times10^{-2}$ at $n_B=1024$; the
phase maximum at the final level is $1.69\times10^{-2}$.  The fitted rate is
$0.400$, compared with $\alpha-\beta=0.31$, and
$\max\operatorname{cond}(\mathbb A_h)=3.198$.  The final classical-spline
error is approximately $120$ times larger.  Figure~\ref{fig:test3-convergence}
shows the convergence for all three tested values of $\beta$.

\subsection{Lateral-data reconstruction at limiting H\"older regularity}

Test~4 uses the contour, jump set, kernel, and cut of Test~2, but places the
H\"older cusp exactly at the second jump, $\theta_c=\theta_2^d=0.83\pi$.
To avoid an artificial antipodal discontinuity, the cusp is localized:
\begin{equation}\label{eq:localized-cusp}
 u_{\mathrm{cusp}}(\theta)
 =\chi(|s|)A_{\pm}
 \left(2\left|\sin\frac{s}{2}\right|\right)^\alpha,
 \qquad
 s=\operatorname{Arg}\bigl(e^{i(\theta-\theta_2^d)}\bigr),
\end{equation}
where $A_-=0.20-0.06i$, $A_+=-0.07+0.16i$, and $\alpha=0.79$.  The cutoff
$\chi$ equals one for $|s|\le0.118\pi$, vanishes for
$|s|\ge0.236\pi$, and is joined by the quintic smoothstep
$1-10x^3+15x^4-6x^5$ on the transition interval.  The localized definition \eqref{eq:localized-cusp} preserves the exact limiting H\"older behavior at the selected jump while removing the spurious discontinuity at the antipodal point.

The samples are taken on a uniformly shifted parameter grid with $N_s=10n_B+1$.  Consequently, the distance to the sixth nearest sample on either adjacent arc is $O(h)$.  The one-sided values are reconstructed from these six nearest samples by the triangular weighted estimator in Algorithm~\ref{alg:lateral}.  Table~\ref{tab:lateral-test4} shows that the
maximum jump error decreases from $3.09\times10^{-2}$ to
$1.82\times10^{-3}$.  The fitted rate is $0.811$, close to the predicted exponent $\alpha=0.79$.  The small difference is consistent with
preasymptotic effects and the smooth cutoff away from the jump.  The corrected full collocation solution also displays decreasing errors; at $n_B=1024$ and $m=4$,
its phase-median errors are $6.29\times10^{-4}$ in the supremum norm and
$8.99\times10^{-3}$ in $PH_{0.45}$.

\begin{table}[t]
\centering
\caption{Test~4: one-sided reconstruction at limiting H\"older regularity.}
\label{tab:lateral-test4}
\small
\begin{tabular}{rcc}
\toprule
$n_B$ & $\max_j|\widehat\delta_j-\delta_j|$ & local rate\\
\midrule
32   & $3.088\times10^{-2}$ & ---\\
64   & $1.744\times10^{-2}$ & $0.824$\\
128  & $9.726\times10^{-3}$ & $0.843$\\
256  & $5.577\times10^{-3}$ & $0.802$\\
512  & $3.211\times10^{-3}$ & $0.796$\\
1024 & $1.824\times10^{-3}$ & $0.816$\\
\bottomrule
\end{tabular}
\end{table}

\begin{figure}[t]
\centering
\includegraphics[width=0.94\textwidth]{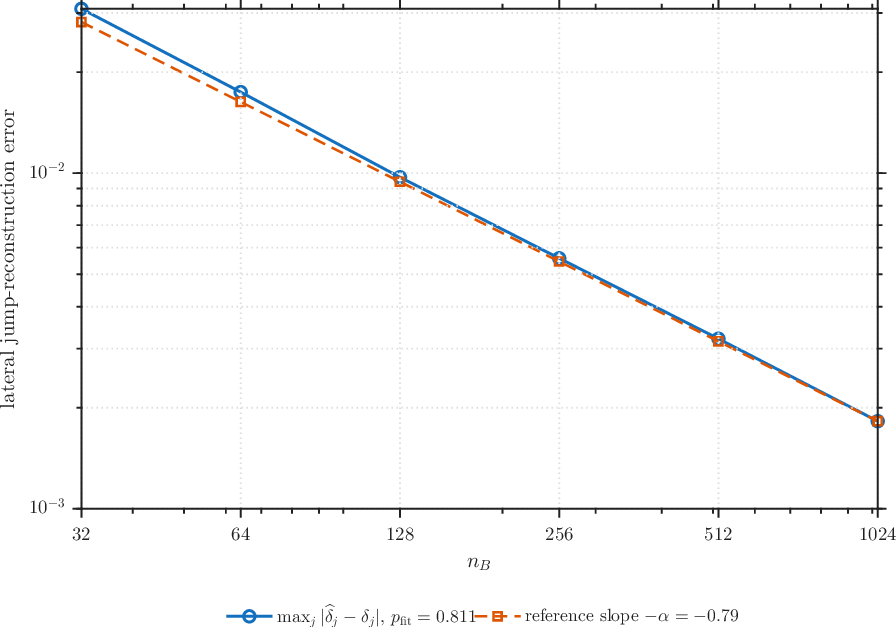}
\caption{Test~4, $m=4$: convergence of the reconstructed lateral data and
jump amplitudes.  The reference slope is $n_B^{-\alpha}$ with $\alpha=0.79$.}
\label{fig:test4-lateral}
\end{figure}

Figure~\ref{fig:test4-lateral} displays the corresponding log--log behavior.

\subsection{Comparison with discontinuity-adapted alternatives}
\label{sec:comparison-results}

At the final target solved dimension $1024$, the geometric means of the raw $PH_\beta$ diagnostics over Tests~1--4 and $m=2,3,4$ are
\[
 9.93\times10^{-3}\quad\text{(BSH)},\qquad
 9.82\times10^{-3}\quad\text{(AGS)},\qquad
 1.22\times10^{-2}\quad\text{(Joe)}.
\]
Thus BSH and AGS have essentially the same aggregate raw accuracy, while
Joe is moderately less accurate in the principal strong-norm diagnostic.  The principal BSH convergence tables report medians over three admissible phases.  The separate comparison module instead uses the middle admissible phase for BSH and closely matched actual solved dimensions for the competing local spaces; this convention accounts for the small differences between the BSH entries below and the phase-median values reported earlier.
The representative order-$4$ values are listed in
Table~\ref{tab:raw-comparison-m4}.  BSH has the smallest error in
Tests~1--3; AGS is better in Test~4, where the limiting H\"older cusp lies
exactly at a physical arc endpoint.

\begin{table}[t]
\centering
\caption{Raw $PH_\beta$ diagnostics for $m=4$ and target solved dimension $1024$.}
\label{tab:raw-comparison-m4}
\small
\begin{tabular}{ccccc}
\toprule
Test & $\beta$ & BSH & AGS & Joe\\
\midrule
1 & $0.55$ & $1.895{\times}10^{-2}$ & $2.122{\times}10^{-2}$
  & $3.048{\times}10^{-2}$\\
2 & $0.45$ & $1.275{\times}10^{-2}$ & $1.973{\times}10^{-2}$
  & $3.150{\times}10^{-2}$\\
3 & $0.50$ & $8.023{\times}10^{-3}$ & $1.055{\times}10^{-2}$
  & $1.626{\times}10^{-2}$\\
4 & $0.45$ & $6.737{\times}10^{-3}$ & $1.736{\times}10^{-3}$
  & $2.367{\times}10^{-3}$\\
\bottomrule
\end{tabular}
\end{table}

One iteration reduces the error by several orders of magnitude.  The
geometric means of the final $PH_\beta$ errors are
\[
 3.02\times10^{-8}\ \text{(BSH)},\qquad
 1.69\times10^{-8}\ \text{(AGS)},\qquad
 7.65\times10^{-9}\ \text{(Joe)},\qquad
 8.39\times10^{-9}\ \text{(Nystr\"om)}.
\]
For BSH this is an aggregate reduction by approximately
$3.3\times10^5$, consistent with the smoothing mechanism in
Theorem~\ref{thm:bsh-iterated}.  Table~\ref{tab:iterated-comparison-m4}
also shows that no iterated method is uniformly best: BSH is best in
Test~2, is close to Joe--Nystr\"om in Test~1, while AGS and the
Gauss-panel methods have smaller constants in the cusp-aligned Test~4.

\begin{table}[t]
\centering
\caption{Iterated $PH_\beta$ diagnostics for $m=4$ and target solved dimension $1024$.}
\label{tab:iterated-comparison-m4}
\small
\begin{tabular}{cccccc}
\toprule
Test & $\beta$ & BSH & AGS & Joe & Nystr\"om\\
\midrule
1 & $0.55$ & $3.497{\times}10^{-8}$ & $6.392{\times}10^{-8}$
  & $3.390{\times}10^{-8}$ & $3.390{\times}10^{-8}$\\
2 & $0.45$ & $2.219{\times}10^{-8}$ & $2.453{\times}10^{-8}$
  & $3.950{\times}10^{-8}$ & $3.950{\times}10^{-8}$\\
3 & $0.50$ & $7.249{\times}10^{-8}$ & $2.522{\times}10^{-8}$
  & $7.106{\times}10^{-8}$ & $6.745{\times}10^{-8}$\\
4 & $0.45$ & $1.220{\times}10^{-8}$ & $5.699{\times}10^{-10}$
  & $7.268{\times}10^{-10}$ & $7.268{\times}10^{-10}$\\
\bottomrule
\end{tabular}
\end{table}

The representations differ more clearly in their jump behavior.  BSH and
AGS reproduce the prescribed jumps to numerical precision in the
constant-leading-coefficient tests.  Raw Joe collocation may jump at every
panel interface; over the twelve final-level configurations, its largest
physical-jump error is $1.04\times10^{-3}$ and its largest artificial
interface jump is $2.86\times10^{-3}$.  Hence its reported raw
$PH_\beta$ value is a discrete diagnostic whenever such additional jumps
are present.  All tested systems remain moderately conditioned through
$n_B=1024$, with no systematic growth under refinement.

The comparison supports the intended positioning of the proposed method:
BSH is not uniformly superior to every adapted alternative, but it combines
competitive raw accuracy with one global periodic spline component, exact
representation of the prescribed physical jumps, and no mesh-dependent
artificial discontinuities.  The complete comparison data and figures are
available in the Zenodo archive cited above.

\subsection{Quadrature check}

The experiments use $q=10$, well above the minimum
$q\ge\lceil m/2\rceil$ in Proposition~\ref{prop:explicit-quadrature}.
For Tests~1, 2, and 4, the difference between the order-$10$ and order-$18$
panel rules is close to floating-point precision; even the isolated
order-$4$ calculation gives defects of order $10^{-14}$--$10^{-15}$.
Test~3 is more informative because $c^{-1}$ is only piecewise H\"older.
Its phase-median solution defect in $PH_{0.50}$ decreases from
$3.01\times10^{-7}$ at $n_B=32$ to $1.14\times10^{-10}$ at $n_B=1024$,
which is negligible relative to the approximation error.  A separate
fixed-mesh calculation with $n_B=256$, $q=4,6,8,10,14$, and reference order
$30$ gives defects between $4.7\times10^{-9}$ and $4.6\times10^{-8}$.
These computations are consistent with the scale required by Theorem~\ref{thm:quadrature}; they are not used to infer a rate with respect
to $q$.

\section{Conclusions}\label{sec:conclusions}

We proposed a periodic B-spline--Heaviside collocation method for regular
Fredholm equations on a closed contour with piecewise H\"older data.  The
central construction combines one global periodic B-spline space with
cut-compensated generators $G_j=H_j-R$, which carry the prescribed physical
jumps without producing a spurious discontinuity at the parameter cut.
Because the regular integral term is globally continuous, the jumps of the
solution equal those of the right-hand side.  The enriched scheme is
therefore triangular: the jump component is reconstructed first, and one
$n_B\times n_B$ collocation system determines the continuous component.

For
$
        0<\beta<\alpha<\mu\le1,
$
the exact method is uniformly stable and converges in $PH_\beta$ at the
rate $O(h^{\alpha-\beta})$.  The fully discrete analysis provides an explicit quadrature estimate:
if the parametrized kernel density is $\rho$-H\"older in the integration
variable with values in $H_\mu$, then a fixed panelwise Gauss rule with
$q\ge\lceil m/2\rceil$ produces an $O(h^\rho)$ perturbation.  Hence
\[
        \norm{\varphi-\varphi_h^q}_{PH_\beta}
        \le
        C\bigl(h^{\alpha-\beta}+h^\rho\bigr),
\]
and the exact-collocation rate is retained whenever
$\rho\ge\alpha-\beta$.  The one-step iterated approximation satisfies the
stronger estimate
\[
        \norm{\varphi-\varphi_h^{\mathrm{it}}}_{H_\mu}
        \le
        C h^\alpha\norm{f}_{PH_\alpha},
\]
without additional assumptions on the data beyond those used for exact
collocation.

The sampled lateral-data reconstruction and the transformation
$w=c\varphi$ extend the triangular construction to practically relevant
input formats and to nonvanishing piecewise H\"older leading coefficients.
The numerical experiments are consistent with the theoretical rates,
moderate conditioning, and the predicted smoothing effect of iteration.
The comparison with AGS, Joe, and Nystr\"om alternatives shows that BSH has
competitive raw accuracy.  Its principal structural advantage is the
simultaneous preservation of a single global periodic spline component,
exact representation of the prescribed jump amplitudes in the constant-leading-coefficient formulation (and in the transformed unknown for the variable-coefficient problem), and absence of mesh-dependent artificial jumps.  Complete secondary tables and figures are retained in
the reproducibility archive so that the main article remains focused on the
B-spline--Heaviside construction and its analysis.

\end{document}